\documentclass[12pt]{article}
\usepackage{amsthm,amsmath,amsfonts,epsfig,amssymb,latexsym}
\usepackage[T1]{fontenc}
\usepackage{graphicx}
\usepackage{enumerate}
\usepackage{xy,xypic}
\usepackage{graphics}
\usepackage{footnote}
\usepackage{scalefnt}
\usepackage{tikz}
\usepackage{color}
\usetikzlibrary{shapes,arrows,fit,calc,positioning}
\usetikzlibrary{matrix}
\xyoption{all}
\title{On D\'{e}vissage for Witt groups} 
\author{
 Satya Mandal \footnote{Partially supported by a General Research Grant from KU}
 ~and  Sarang Sane\\ %$^\dagger$\\ 
{\small University of Kansas, Lawrence KS 66045}\\
{\small {\it emails: mandal@math.ku.edu, 
 ssane@math.ku.edu}}\\
{\small Keywords : d\'{e}vissage, Witt groups, derived categories}\\
{\tiny Mathematics Subject Classification[2000]: Primary:11E81, 18E30, 19G12 ;
 Secondary: 13C10, 13D}\\
}
\begin{document}
\renewcommand{\baselinestretch}{1.255}
\setlength{\parskip}{1ex plus0.5ex}
\date{24 November, 2013}
\newcommand{\iso}{\stackrel{\sim}{\longrightarrow}}
\newcommand{\sur}{\twoheadrightarrow}
\newcommand{\bp}{\begin{proposition}}
\newcommand{\ep}{\end{proposition}}
\newcommand{\bl}{\begin{lemma}}
\newcommand{\el}{\end{lemma}}
\newcommand{\bt}{\begin{theorem}}
\newcommand{\et}{\end{theorem}}
\newcommand{\eop}{\hfill \square}
\newcommand{\pf}{\noindent{\bf Proof.~}}
\newcommand{\PD}{\text{proj} \dim}
\newcommand{\lra}{\longrightarrow}
\newcommand{\Lra}{\Longrightarrow}
\newcommand{\pic}{The proof is complete.}
\newcommand{\Dia}{\diagram}
%%%%%%%%%%%%%%%%SSSSSSSSSSSSSSSSSSSSSSSSSSSSSSSSSSSSS I've added these
\def\spec#1{\mathrm{Spec}(#1)}
\def\m{\mathfrak {m}}
\def\A{\mathcal {A}}
\def\B{\mathcal {B}}
\def\P{\mathcal {P}}
\def\C{\mathcal {C}}
\def\D{\mathcal {D}}
\def\E{\mathcal {E}}
\newcommand{\smallcirc}[1]{\scalebox{#1}{$\circ$}}
\def\Z{\mathbb {Z}}
%%%%%%%%%%%%%%%%SSSSSSSSSSSSSSSSSSSSSSSSSSSSSSSSSSSSS Done
\newcommand{\bE}{\begin{enumerate}}
\newcommand{\eE}{\end{enumerate}}
\newcommand{\bT}{\begin{theorem}} 
\newcommand{\eT}{\end{theorem}}

\newtheorem{theorem}{Theorem}[section]
\newtheorem{proposition}[theorem]{Proposition}
\newtheorem{lemma}[theorem]{Lemma}
\newtheorem{corollary}[theorem]{Corollary}
\newtheorem{construction}[theorem]{Construction}
\newtheorem{notations}[theorem]{Notations}
\newtheorem{question}[theorem]{Question}
\newtheorem{example}[theorem]{Example}
\newtheorem*{theorem*}{Theorem}
\theoremstyle{remark}
\newtheorem{rem}[theorem]{Remark}
\theoremstyle{definition}
\newtheorem{defin}[theorem]{Definition}
\maketitle

%%%%%%%%%%%%%%%%%%%%%%%%%%%%%%%%%%%%%%%%%%%%%%%%%%%%%%%%%%%%%%%%%%
%%%%%%%%%%%%%%%%%%%%%%%%%%%%%%%%%%%%%%%%%%%%%%%%%%%%%%%%%%%%%%%%%%
\section{Introduction} 
%%%%%%%%%%%%%%%%%%%%%%%%%%%%%%%%%%%%%%%%%%%%%%%%%%%%%%%%%%%%%%%%%%

%%%SSSSSSSSSSSSSSSSSSSSSSSSSSSSSSSSSSSSSSSSS
%%%%%%%%%%%%%%%%%%%%%%%%%%%%%%%%%%%%%%%%%%%%%%%%%%%%%%%%%%%%%%%%%%Introduction is rewritten
In this article, we develop a version of the d\'{e}vissage theorem \cite[Theorem 6.1]{BW}
for Witt groups of Cohen-Macaulay rings. To introduce this theorem, let $A$ be a commutative
Noetherian domain of dimension $d$ with $2$ invertible and $K$ be its quotient field. It is a
classical question (known as purity) as to when the map $W(A) \rightarrow W(K)$ is injective.
Purity is a conjecture when $A$ is a regular local ring and is affirmatively settled in  (\cite{CS, OP, OPSS}). %\textcolor{blue}{many cases}.

%%Citations needed.

In general, we can extend the above map to the right for any regular scheme by considering the
Gersten-Witt complex. Let $X=\spec{A}$ be of dimension $d$ with $2$ invertible and $X^{(n)}$ denote the points of codimension $n$. A Gersten-Witt complex $$ 0 \rightarrow W(A) \rightarrow 
\bigoplus_{x \in X^{(0)}} W(k(x)) \rightarrow \bigoplus_{x \in X^{(1)}} W(k(x)) \rightarrow
\ldots \rightarrow \bigoplus_{x \in X^{(d)}} W(k(x)) \rightarrow 0$$ was first constructed
by Pardon \cite{Pa1}. He further conjectured the exactness of this sequence for regular local rings and affirmatively settled it in 
 many cases.

Subsequently, with the introduction of triangular Witt groups by Balmer \cite{DWG, TWGI,TWGII},
Witt groups could be viewed as a cohomology theory. Using this, another Gersten-Witt complex could be defined (though both complexes look similar, it seems unproven that the differentials match) \cite{BW}, similar to the one in K-theory. Once again it is an open question as to when the complex is exact. In particular, it is conjectured to be so when $X=\spec{A}$ where $A$ is a regular, local ring and this is affirmatively settled in  (\cite{B4, Pa1, Pa2, B4, BGPW}).
%\textcolor{blue}{many cases}. 
\cite{DWG, TWGI,TWGII, BW} form the basic backbone 
of the methodologies  in this article and we would like to remark that the interested reader would be well advised to take a look at them. {\it For any unexplained 
notations and definitions in the introduction, please refer to (\ref{notaCate})}.

The key result which allows one to move from derived Witt groups to the Gersten-Witt complex is
d\'{e}vissage~\cite[Section 6]{BW} which states that for a regular, local ring $(A,\m ,k)$, of
dimension $d$, we have $W^n(D^b_{fl}(A)) \cong W(k)~ \text{if}~ n \equiv d ~\text{mod}~ 4~ \text{and}~W^n(D^b_{fl}(A)) \cong 0 $ otherwise. We describe below the generalized form of this theorem for a Cohen-Macaulay ring $A$ with $\dim A_{\m}=d$ for all maximal ideals $\m$ and $2$ invertible.

Suppose $A$ is a Cohen-Macaulay ring with $\dim A=d$. Since there are modules of infinite
projective dimension over such a ring and of finite length, the usual duality $Ext^d(~\underline{~~}~,~A)$
does not work well. The options are either to change the duality (w.r.t. the canonical module) but then use the category of all modules (coherent Witt groups) or impose finite projective dimension homology conditions on the complexes. The first path is taken and deeply studied
in Pardon~\cite{Pa2, Pa3} and more recently in Gille~\cite{G1,G2}, where they also establish
a Gersten-Witt complex of coherent Witt groups.

In this article, we take the second approach. Let ${\mathcal M}FPD(A)$ be the category of
finitely generated $A-$modules with finite projective dimension, ${\mathcal M}FL(A)$ be
the category of finitely generated $A-$modules with finite length, and ${\mathcal A}={\mathcal M}FPD_{fl}(A)$ be the full subcategory of finitely generated $A-$modules with finite projective dimension and finite length (the "intersection"). Note then that $\A$ is an exact category
and has a natural duality given by $M \mapsto Ext^d_A(M,A)$ and so we can consider the Witt
group $W(\A)$. By Balmer~\cite{TWGI}, we already know that $W({\A}) \cong W^d(D^b({\A}))$.

We consider the category $D^b_{\A}(\A)$ with homologies in $\A$. Based on the fact that $\A$
actually has the $2$-out-of-$3$ property for objects, we prove that the duality actually
restricts to $D^b_{\A}(\A) $ and this allows a suitable definition of Witt groups 
$W^i(D^b_{\A}(\A))$. Once defined, we prove that the above isomorphism actually factors
through isomorphisms $W(A) \iso W^d(D^b_{\A} ({\A})) \iso W^d(D^b ({\A}))$.

We now consider the category $D^b_{\A}({\P}(A))$. Similar to $D^b_{\A}(\A)$, we establish
that this category is stable  
 under duality and define the Witt groups $W^i(D^b_{\A}(\A))$.
Having done so, we prove our version of d\'{e}vissage (\ref{IsoOfWzeta}, \ref{shiftFinal}) :
\begin{align*}
W({\A}) & \iso W^d(D^b_{\A}({\P}(A))) \\
W^{-}({\A}) & \iso W^{d+2}(D^b_{\A}({\P}(A))) \\
W^{d+1}(D^b_{\A}({\P}(A))) & \cong W^{d-1}(D^b_{\A}({\P}(A))) \cong 0.
\end{align*}

Note that when $A$ is regular, this is exactly the same theorem as that in \cite{BW}.
Further, in the Gorenstein case, there is a natural commuting square 
involving the above d\'{e}vissage statement 
and the d\'{e}vissage statement in \cite{G1}, %\textcolor{magenta}{
for coherent Witt groups. As of now,
we do not know if these sets of groups coincide %\textcolor{magenta}{
or not,
which  would also be a subject of future investigations. 
However, we point out that there are more forms and neutral spaces in the realm of coherent Witt theory. Since we prove the theorem without
regularity assumptions or the existence of a natural dualizing complex, we do not have
access to the equivalences of the derived categories with duality as in \cite{BW} or \cite{G1}
(in particular we cannot use the powerful lemma of Keller \cite[\textsection 1.5, Lemma and Example(b)]{K}). Our proof is thus necessarily more elementary and na\"{i}ve than the one in \cite{BW}. One of the key ingredients in the proof is the construction of a special sublagrangian
(\ref{leftSubLag}) for symmetric forms in $(D^b_{\A}({\P}(A)))$. 
%\textcolor{blue} 
%{All our results are heavily
%dependent on the beautiful papers of Balmer \cite{DWG, TWGI,TWGII} and we again remind
%the reader to take a look at them.
%} 

%\textcolor{purple}
In deed, the methods in this paper have  much wider applications. 
This is the first of a series of articles (\cite{MS2, MS3})
dedicated to apply the sublagrangian theorem
of Balmer (\cite{TWGII}), or otherwise,  
to singular varieties and provide further insight into nonsingular varieties.
Our interest in these studies  stems from the introduction of the
Chow-Witt groups $\widetilde{CH}^r(A)$ for $0\leq r\leq d$, due to Barge and Morel
\cite{BM} and developed by Fasel \cite{F1}, as 
the obstruction groups 
for splitting of projective modules, 
%\textcolor{magenta}{
which works best for nonsingular 
varieties.
This also serves as a motivation for our 
interest in
maintaining the category of projective
modules in our statement of d\'{e}vissage.
%and the reluctance to move to coherent Witt groups. 
Jean Fasel informs that, for singular varieties $X$,
Chow-Witt groups $\widetilde{CH}^r(X)$ and obstruction classes
can be  defined in the same manner, using coherent Witt groups. 
However, it is not known whether vanishing of the obstruction classes
would lead to splitting.
We feel that, for the purpose of developing
an obstruction theory for singular varieties,
 it would be more natural to consider 
 some analogue of the derived Witt groups of the category of projective modules.
This approach would be of its own 
independent interest for both computations and applications.
% }
% 
%\textcolor{blue}{We feel that, for the purpose of developing
%an obstruction theory with this property  for singular varieties,
% it is more natural to use some analogue of the derived Witt groups. 
%}
%
%While it would be possible to give a definition of the Chow-Witt groups, 
%using the Coherent Witt groups, we are not convinced that such an approach 
%may lead to a meaningful obstruction theory for singular varieties.}
%\textcolor{blue}{
%Since there is a parallel theory of Euler class groups, which is defined over
%a Cohen-Macaulay ring, the eventual goal
%would be to define Chow-Witt groups also over
%any Cohen-Macaulay ring.
%}
%\textcolor{magenta}{
The  results in this article also give 
rise to the possibility of a Gersten-Witt complex for these Witt groups,
%\textcolor{magenta}{
parallel to that of coherent Witt groups.

A word about the layout of the article : in section~\ref{Background}, we establish the
basic definitions and a key result on projective dimensions. In the section~\ref{Duality},
we establish the important theorem that the categories $D^b_{\A}({\P}(A))$ and $D^b_{\A}({\A})$
are closed under duality, and more specifically how the homologies of the dual look like. Once
this is established, in section~\ref{WittGroups}, we define the Witt groups of the above categories and expectedly, they are $4$-periodic, i.e. $W^n(D^b_{\A}({\P}(A))) \iso W^{n+4}(D^b_{\A}({\P}(A)))$. Finally, in sections \ref{WDIAGRAM} and \ref{shiftedWSec},we prove our main theorems about d\'{e}vissage.

\noindent{\bf Acknowledgement:} {\it We would like to thank Sankar P. Dutta
for many helpful discussions, and particularly, for example (\ref{Dutta}). 
We are also very thankful to Jean Fasel for many helpful discussions and guidance. Thanks
are also due to Manuel Ojanguren for some helpful corrections and comments. The second named author is thankful to the University of Kansas and the Robert D.Adams Trust.}

\section{Basic Notations and Preliminaries}\label{Background}  

\begin{notations}\label{notaCate}{\rm 
Throughout this article, 
{\bf $A$ will denote a Cohen-Macaulay ring with $\dim A_m=d \geq 2$,
for all maximal ideals $m$ of $A$. Further, $2$ is always invertible in $A$.}
We set up the notations :
\begin{enumerate}
\item ${\mathcal M}(A)$ : category of finitely generated $A-$modules.
\item ${\mathcal M}FPD(A)$ : full subcategory of finitely generated $A-$modules 
with finite projective dimension.
\item ${\mathcal M}FL(A)$ : category of finitely generated $A-$modules with finite length.
\item  $\A = {\mathcal M}FPD_{fl}(A)$ : category of finitely generated $A-$modules with
 finite length and finite projective dimension.
\item ${\mathcal P}(A)$ : category of finitely generated projective $A-$modules.
\item For any exact category $\C$, $Ch^b(\C)$ is the category of bounded chain complexes
with objects in $\C$, and $D^b(\C)$ is its derived category.
\item For any two exact categories $\C, \D$ in an ambient abelian category $\C'$,
$Ch^b_{\D}(\C)$ is the full subcategory of $Ch^b(\C)$ consisting of complexes
with homologies in $\D$. $D^b_{\D}(\C)$ is its derived category, which is also the
full subcategory of $D^b(\C)$ consisting of objects from $Ch^b_{\D}(\C)$.
\item ${\mathcal R}$ : full subcategory of $D^b_{{\A}}({\mathcal P}(A))$ consisting of objects
$P_{\bullet}$ such that $P_i=0$ for $i > d, i < 0$ and $H_i(P_{\bullet})=0$ for all $i\neq 0$
and $H_0(P_{\bullet})\in {\mathcal A}$.
\item  For objects $M$ in ${\mathcal A}$,
let $M^{\vee}=Ext_A^d(M, A)$ and $\tilde{\varpi}: M \iso M^{\vee\vee}$
be the identification by double ext. (but we make this more precise in 
diagram (\eqref{diag1}) and the explanation of $\iota$.)
\item We will denote complexes $P_{\bullet}$ by : \\
$
\diagram
\cdots 0 \ar[r] & P_m \ar[r]^{\partial_m}
&P_{m-1} \ar[r] & \cdots \cdots \ar[r]  &P_n \ar[r] & 0 \cdots
\enddiagram 
$
\item A non-zero complex $P_{\bullet}$ is said to be supported on $[m,n]$ if
$P_i=0$ for all $i<n$ and $i>m$.
\item For a complex  $P_{\bullet}$ of projective $A-$modules
$P_{\bullet}^*$ will denote the usual dual induced
by $Hom(*,A)$ and $\varpi:P_{\bullet}\iso P_{\bullet}^{**}$ will denote the
identification by evaluation. Note that the degree $r-$component of the dual 
$P_{\bullet}^*$ is $(P_{-r})^*$.
\item The {\bf length} of a non-zero complex $P_{\bullet}$ is defined as
$\ell(P_{\bullet})=u-l$ where $P_u\neq 0, P_l\neq 0$ and $P_i=0$ for all
$i<l$ and $i>u$.
\item Let $B_r=B_r(P_{\bullet}):= \partial_{r+1}(P_{r+1}) \subseteq P_r$
denote the module of  $r-$boundaries and $Z_r=Z_r(P_{\bullet}):= \ker(\partial_{r}) \subseteq P_r$
denote the module of $r-$cycles (or the $r^{th}$ syzygy).
\item The $r^{th}-$homology of $P_{\bullet}$ will be denoted by $H_r=H_r(P_{\bullet}):= \frac{Z_r}{B_r}$. So, the $r^{th}-$homology of the dual is
$H_r(P_{\bullet}^*)=\frac{\ker\left(\partial_{-(r-1)}^*\right)}{image(\partial_{-r}^*)}$. 
\item A full exact subcategory $\C$ of an abelian category $\D$ is said to have the
$2$-out-of-$3$ property if for every short exact sequence in $\D$, whenever two of the
objects are objects of $\C$, then so is the third.
\end{enumerate}
}
\end{notations}

\begin{rem}\label{rmk1}
The subcategory ${\mathcal M}FPD(A)$ and $\A$ are both exact subcategories and in fact
have the $2$-out-of-$3$ property. The category ${\mathcal R}$ is also an exact category.
Although it is a subcategory of $D^b_{{\A}}({\mathcal P}(A))$, it has no translation
and is actually naturally equivalent to the category $\A$. The natural functor 
$\eta : {\mathcal R} \rightarrow \A$ is given by sending a complex $Q_{\bullet}$ to
$H_0(Q_{\bullet})$. The inverse functor $\iota$ is given by associating to objects
$M \in {\mathcal A}$ a projective resolution of length $d$.

Note further that when $A$ is not regular, the categories $D^b_{{\A}}({\mathcal P}(A))$
and $D^b_{{\mathcal M}FPD(A)}({\mathcal P}(A))$ are not closed under the cone operation 
as the following example demonstrates. 
\begin{example}\label{failExact}{\rm
Let $(A,\m)$ be a non-regular Cohen-Macauly ring with $\dim A=d$, such that 
$\m=(f_1, f_2, \ldots, f_d, z)$. We can assume, using prime avoidance,
that $f_1, f_2, \ldots, f_d$ is a regular sequence. Let
$U_{\bullet}=Kos_{\bullet}(f_1, f_2, \ldots, f_d)$ be the Koszul complex.
Since the only nonzero homology of $U_{\bullet}$ is 
$H_0(U_{\bullet})=\frac{A}{(f_1, f_2, \ldots, f_d)} \in {\mathcal A}$,
$U_{\bullet}$ and all its translates are objects of both the categories above.
Let $C(z)$ denote the cone of the the chain complex map $z:U_{\bullet} \lra U_{\bullet}$.
From the long exact homology sequence corresponding to the short exact sequence of
chain complexes
$$
\diagram 
0 \ar[r] & U_{\bullet} \ar[r] & C(z) \ar[r] & U_{\bullet}[1] \ar[r] & 0
\enddiagram 
$$
it follows that
$$H_0(Cone(z)) \cong \text{coker}(\frac{A}{(f_1, f_2, \ldots, f_d)} 
\stackrel{\cdot z}{\rightarrow}
\frac{A}{(f_1, f_2, \ldots, f_d)}) \cong \frac{A}{\m} \notin {\mathcal A}.$$
So, $C(z)$ is not an object in $D^b_{{\mathcal A}}({\mathcal P}(A)$.
}
 \end{example}
\end{rem}

Next, we ask if $Ch^b_{{\mathcal M}FPD}({\mathcal P}(A))$ is closed under duality. We thank
Sankar Dutta for providing the following example :

\begin{example}[Dutta]\label{Dutta}{\rm
Let $(A,\m,k)$ be any non-regular Cohen-Macaulay local ring,
with $\dim A=d$. Let 
$$
\diagram
\cdots \ar[r] & P_d \ar[r]^{\partial_d}
&P_{d-1} \ar[r] & \cdots \ar[r] & P_0 \ar[r] &k \ar[r] &0
\enddiagram 
$$
be a projective resolution of $k$. Let $^*$ denote $Hom(-,A)$ and    
$M=cokernel(\partial_d^*)$. Since $Ext^r(k, A)=0$ for all $0\leq r<d$, the 
sequence 
$$
\diagram
0 \ar[r] & P_0^* \ar[r] & \cdots \ar[r]
& P_{d-1}^*\ar[r] & P_d^* \ar[r] &M \ar[r] &0
\enddiagram 
$$
is a projective resolution of $M$. Dualizing this sequence, it follows that
$Ext_A^d(M,A) \cong k$, which does not have finite projective dimension.
In particular,  $Ch^b_{{\mathcal M}FPD}({\mathcal P}(A))$
is not closed under duality.
}
\end{example}

Note however that in the above example, $M$ does  not have finite length. Indeed, we
will prove in section~\ref{Duality} that the category $Ch^b_{{\A}} ({\P} (A))$ is closed
under duality.

We mention a few standard results for the sake of completeness.
\begin{lemma}\label{infExt}
Let  $(A, m)$ be a Cohen-Macaulay local ring with $\dim A =d.$  
Let $M \in {\mathcal M}FL(A)$. Then Further, $Ext^i(M,A)=0$ for all $i < d$.
and $Ext^d(M, A) \neq 0$ is also in ${\mathcal M}FL(A)$.
Note further that if $M \in \A$, then so is $Ext^d(M, A)$.
\end{lemma}
\begin{lemma}\label{ExtDual} 
Let  $A$ be a Cohen-Macaulay  ring with $\dim A =d.$
Let $M\in {\mathcal A}.$
There is a natural isomorphism
$$
\varpi: M \iso M^{\vee\vee}.  
$$
\end{lemma}

\begin{corollary}\label{corPflCat} 
$({\mathcal A}, ^{\vee}, \pm \tilde{\varpi})$ are  exact categories with duality.
\end{corollary}

%%%%%%%%%%%%%%%%%%%%%%%%%%%%%%%%%%%%%%%%%%%%%%%%%%%%%%%% The propn. and the proof below needs serious rewriting. In particular, $A$ has to be CM with all maximal ideals of ht. $d$.

\begin{proposition}\label{PDofBrZr} 
Let $A$ be Cohen-Macaulay with $d=\dim A_{\m} \geq 2$ for all maximal ideals $\m$.
Let $P_{\bullet}$  be a complex in $Ch^b({\mathcal P}(A))$. Assume that all the homologies $H_r:=H_r(P_{\bullet}) \in {\mathcal A}$. Then we claim :
\begin{enumerate}
\item \label{allFpd} The modules $B_r$ and $Z_r$ have finite projective dimension $~\forall r$.
In that case, $proj\dim(Z_r)= proj\dim(B_{r-1})-1$.
\item \label{ZrAtmost} For all $r \in {\mathbb Z}$, we have $proj\dim B_r \leq d-1$
and so $proj\dim Z_r \leq d-2$.
\item \label{HrNot0} If $H_r\neq 0$ then $proj\dim B_r=d-1.$
\end{enumerate} 
\end{proposition}
\noindent{\bf Proof.} Note that $P_{\bullet}$ is a bounded complex and so let it be
supported on $[m,n]$. Then $Z_n = P_n$ and so has finite projective dimension.
Since there are short exact sequences :
$$0 \rightarrow B_r \rightarrow Z_r \rightarrow H_r \rightarrow 0 ~~~~~~~~~~
0 \rightarrow Z_r \rightarrow P_r \rightarrow B_{r-1} \rightarrow 0,$$
it is clear that if $Z_r$ is of finite projective dimension, then so are
$B_r$ and $Z_{r+1}$, and hence the proof follows by induction. The second part of (\ref{allFpd})
also follows from the above exact sequence.

Since $B_r$ is torsion free, it has depth at least $1$ and hence, by the Auslander-Buchsbaum
theorem, $\PD(B_r) \leq d-1$. So, $\PD(Z_r) \leq d-2$. So, (\ref{ZrAtmost}) is established. 

Now, assume $H_r \neq 0$. Choose a maximal ideal $\m$ in the support of $H_r$. Then, consider
the localized short exact sequence 
$$0 \rightarrow (B_r)_{\m} \rightarrow (Z_r)_{\m} \rightarrow (H_r)_{\m} \rightarrow 0.$$
Then, we get a long exact sequence of $Tor_{A_{\m}}(\_,A/\m)$, which gives us that $$Tor^d_{A_{\m}}((H_r)_{\m} ,A/\m) \cong Tor^{d-1}_{A_{\m}}((B_r)_{\m} ,A/\m) \quad
 \text{and} \quad Tor^d_{A_{\m}}((B_r)_{\m} ,A/\m) \cong 0,$$
since we have already proved that $proj\dim_A Z_r \leq d-2$.
Thus we obtain $proj\dim_{A_{\m}} (B_r)_{\m} = d - 1$. Since we know that
$\PD_A(B_r) \leq d-1$, this implies $\PD_A(B_r) = d-1$. This establishes (\ref{HrNot0}).
$\eop$ 

%%%%%%%%%%%%%%%%%%%%%%%%%%%%%%%%%%%%%%%%%%%%%%%%%%%%%%%%%%%%%%%%%%%%%%%%%%%%%%%%%%%%%%%%%%

\vspace{5mm}
The complexes in $Ch^b({\mathcal P}(A))$ with finite length homologies have at least $d$
nonzero components at the left where the homology is $0$. This proposition plays
a key role in sections \ref{WDIAGRAM} and \ref{shiftedWSec}.
%%%%%%%%%%%%%%%%%%%%%%%%%%%%%%%%%%%%%%%%%%%%%%%%%%%%%%%%%%%%%%%%%%%Needs to be rewritten.
\begin{proposition}\label{GoodAtLeft}
Let $A$ be Cohen-Macaulay with $d=\dim A_{\m}$ for all maximal ideals $\m$.
Let $P_{\bullet}$ be a bounded complex of projective modules,
such that $H_i=0 ~\forall i > n$ and $H_n \neq 0$ is of finite length.
Then $P_i \neq 0, n \leq i \leq n+d$.
\end{proposition}
\pf 
Since $H_i = 0 ~\forall i > n$ and the complex is bounded, we get that $\frac{P_n}{B_n}$ is of
finite projective dimension, since the components with indices $\geq n$ give a resolution.
Now, let $\m$ be a maximal ideal in the support of $H_n(P_{\bullet})$, then 
$(H_n(P_{\bullet}))_{\m} \subseteq \frac{(P_n)_{\m}}{(B_n)_{\m}}$ is of finite length,
and hence $\frac{(P_n)_{\m}}{(B_n)_{\m}}$ has depth $0$. By the Auslander-Buchsbaum theorem,
$\PD_{A_{\m}}(\frac{(P_n)_{\m}}{(B_n)_{\m}}) = d$. But that means $\PD_A(\frac{P_n}{B_n}) = d$.
Hence, the resolution of $\frac{P_n}{B_n}$ given by the components of $P_{\bullet}$ with indices
$\geq n$ must have length at least $d$. Hence, $P_i \neq 0, n \leq i \leq n+d$.
\pic $\eop$ 

\section{Duality}\label{Duality}
As always, {\bf $A$ will denote a Cohen-Macaulay ring with $\dim A_m=d \geq 2$, for all maximal
ideals $m$} and ${\mathcal A}={\mathcal M}FPD_{fl}(A)$. In this section, we prove that the
category $Ch^b_{\A}({\P}(A))$ is closed under duality and give a precise description of
the homologies of the dual.

\begin{theorem}\label{TheEasyProof}
Suppose $P_{\bullet}$ is a complex in $Ch^b({\mathcal P}(A))$ with 
homologies in ${\mathcal A}$. Then we have :
\begin{align}
Ext^i(Z_{r}, A) & \cong
\left\{ 
\begin{array}{ll}
Ext^d(H_{r+i-(d-1)},A)& %d-1-r 
1\leq i \leq d-2 \\
0 & for~ i\geq d-1  
\end{array} 
\right. \\
Ext^i(B_{r}, A) & \cong
\left\{ 
\begin{array}{ll}
Ext^d(H_{r+i-(d-1)},A)& %(d-1)-r
1\leq i \leq d-1 \\
0 & for~ i\geq d  
\end{array} \right. \\
Ext^i\left(\frac{P_{r}}{B_{r}}, A\right) & \cong
\left\{ 
\begin{array}{ll}
Ext^d(H_{r+i-d},A)& %d-r
1\leq i \leq d \\
0 & i\geq d  
\end{array} 
\right.
\end{align}
\end{theorem}
\pf
Since $P_i=0~ \forall i \ll 0$, 
the theorem is true for $r\ll 0$. So, we assume that 
the theorem 
is true for $r-1$ and prove it for $r.$

Corresponding to the short exact sequence
$0 \rightarrow Z_{r} \rightarrow P_{r} \rightarrow B_{r-1}\rightarrow 0$,
we get a long exact Ext-sequence which yields
\begin{align}
0 \rightarrow Ext^0(B_{r-1}, A)  \rightarrow P_{r}^* \rightarrow Ext^0(Z_r,A) 
\rightarrow Ext^1(B_{r-1}, A)  \rightarrow 0
\end{align}
and for $i\geq 1$ we have $Ext^i(Z_r,A) \cong Ext^{i+1}(B_{r-1}, A)$. 
Thus, the induction hypothesis yields that for $i\geq 1$,
$$
Ext^i(Z_r, A)= Ext^{i+1}(B_{r-1}, A)=  \left\{ 
\begin{array}{ll}
Ext^d(H_{r+i-(d-1)},A)& %d-1-r 
1\leq i \leq d-2 \\
0 & for~ i\geq d-1  
\end{array} 
\right. \\
$$
So, equation (1) is established.

Consider the long exact Ext-sequence corresponding to the short exact sequence
$0 \rightarrow B_{r} \rightarrow Z_{r} \rightarrow H_{r}\rightarrow 0$.
By (\ref{infExt}), $Ext^i(H_r, A)=0$ for all $i\neq d$ and since 
$Ext^i(Z_r, A)=0~~ \forall i > d-2$ from equation (1), it follows that

\begin{align*}
Ext^i(B_{r}, A) ~~ & \cong ~~
\begin{cases}
Ext^i(Z_r,A),~~ 0\leq i \leq d-2 \\
Ext^d(H_r,A),~~  i=d-1\\
0, ~~~~~~~~~~~~~~~~~~~~~i \geq d
\end{cases} \\
& \cong ~~
\begin{cases}
Ext^0(Z_r,A),~~~~~~~~~ i=0 \\
Ext^d(H_{r+i-(d-1)},A),~~ 1\leq i \leq d-1 \\
0, ~~~~~~~~~~~~~~~~~~~~~~~~~~~~ i\geq d
\end{cases}
\end{align*}

Now consider the short exact sequence
$0 \rightarrow H_{r} \rightarrow \frac{P_r}{B_r} \rightarrow B_{r-1}  \rightarrow 0.$
Again, $Ext^i(H_r, A)=0 ~\forall~i\neq d$ from (\ref{infExt}) and from equation (2) we
get $Ext^i(B_r, A)=0 ~\forall i > d-1$. So it follows that 
\begin{align*}
Ext^i\left(\frac{P_r}{B_r}, A\right) ~~ & \cong ~~
\begin{cases}
Ext^i\left(B_{r-1}, A\right),~~ 0\leq i \leq d-1 \\
Ext^d(H_r,A),~~~~  i=d\\
0, ~~~~~~~~~~~~~~~~~~~~~i > d
\end{cases} \\
& \cong ~~
\begin{cases}
Ext^i\left(B_{r-1}, A\right),~~~~~ i=0 \\
Ext^d(H_{r+i-d},A),~~ 1\leq i \leq d \\
0, ~~~~~~~~~~~~~~~~~~~~~~~~~~~~ i > d
\end{cases}
\end{align*}
$\eop$

\begin{corollary}
Suppose $P_{\bullet}$ is a complex in $Ch^b({\mathcal P}(A))$ with 
homologies in ${\mathcal A}$. Then, for all $r \in {\mathbb Z}$
$$
Ext^i(B_r,A), \quad Ext^i(Z_r,A), \quad Ext^i\left(\frac{P_r}{B_r},A\right)
$$
are in ${\mathcal A}$ for $i\geq 1$ and are in ${\mathcal M}FPD(A)$ for $i=0$.
\end{corollary}
\pf
By (\ref{infExt}) and the preceding theorem (\ref{TheEasyProof}), for $i \geq 1$, the
statement is clear. For $i =0$, we recall below equation (4) from the preceding proof :
\begin{align*}
0 \rightarrow Ext^0(B_{r-1}, A)  \rightarrow P_{r}^* \rightarrow Ext^0(Z_r,A) 
\rightarrow Ext^1(B_{r-1}, A)  \rightarrow 0.
\end{align*}
and that we also proved  $Ext^0(B_{r}, A) \cong Ext^0(Z_{r}, A)$ and
$Ext^0\left(\frac{P_r}{B_r}, A\right) \cong Ext^0\left(B_{r-1}, A\right)$.
Hence, it is enough to know that $Ext^0\left(B_{r-1}, A\right)$ satisfies the theorem.
Once again induction saves the day!
$\eop$

This allows us to conclude our main theorem of the section.
\begin{theorem} \label{easyHLofDual}
Let $P_{\bullet}$ be a complex as in theorem (\ref{TheEasyProof}). Then,
for $t\in{\mathbb Z}$, we have 
$$H_{-t}(P_{\bullet}^*) \cong Ext^d\left(H_{t-d}(P_{\bullet}) , A\right) \cong 
H_{t-d}(P_{\bullet})^{\vee}. $$
In particular, $H_r(P_{\bullet}^*) \in \A$ and hence,
$Ch^b_{\A}({\P}(A))$ is closed under duality.
\end{theorem} 
\pf
Consider the dual complex :
$\cdots P_{t-1}^* \stackrel{(\partial_{t})^*}{\rightarrow}  P_{t}^* 
\stackrel{(\partial_{t+1})^*}{\rightarrow}  P_{t+1}^* \cdots$.
Note that $(\partial_{t+1})^* : P_{t}^*  \rightarrow P_{t+1}^*$
factors through 
$$P_{t}^*  \rightarrow B_{t}^* \hookrightarrow \left(\frac{P_{t+1}}{B_{t+1}}\right)^*
 \hookrightarrow P_{t+1}^*$$
(recall $d \geq 2$) and hence, $ker((\partial_{t+1})^*) = ker(P_{t}^*  \rightarrow B_{t}^*) =
\left(\frac{P_{t}}{B_{t}}\right)^*$.
Similarly, $ker((\partial_{t})^*) = \left(\frac{P_{t-1}}{B_{t-1}}\right)^*$ and hence
we obtain the exact sequence
$$0\rightarrow \left(\frac{P_{t-1}}{B_{t-1}}\right)^* \rightarrow
\left(P_{t-1}\right)^* \rightarrow
\left(\frac{P_{t}}{B_{t}}\right)^* \rightarrow
 H_{-t}(P_{\bullet}^*) \rightarrow 0.$$

Note also that there is an exact Ext-sequence
$$0\rightarrow \left(\frac{P_{t-1}}{B_{t-1}}\right)^* \rightarrow
\left(P_{t-1}\right)^* \rightarrow (B_{t-1})^* \rightarrow
Ext^1\left(\frac{P_{t-1}}{B_{t-1}}, A \right) \rightarrow 0.$$

But since $d \geq 2$, we have
{\scalefont{0.7}
$$
\diagram
0\ar[r] & \left(\frac{P_{t-1}}{B_{t-1}}\right)^* \ar[r] \ar@{=}[d]
& \left(P_{t-1}\right)^* \ar[r]\ar@{=}[d]
& (B_{t-1})^* \ar[d]^{\wr}\ar[r]
& Ext^1\left(\frac{P_{t-1}}{B_{t-1}} ,A\right)\ar@{-->}[d]^{\wr}\ar[r] & 0\\
0\ar[r] & \left(\frac{P_{t-1}}{B_{t-1}}\right)^* \ar[r] 
& \left(P_{t-1}\right)^* \ar[r]
& \left(\frac{P_{t}}{B_{t}}\right)^* \ar[r]
& H_{-t}(P_{\bullet}^*) \ar[r] & 0
\enddiagram
$$
}
Hence, by (\ref{TheEasyProof}), we get that
$H_{-t}(P_{\bullet}^*) \cong Ext^d\left(H_{t-d}(P_{\bullet}) , A\right)$.
The rest follows from (\ref{infExt}).
$\eop$

\begin{rem}\label{rmk1.5nat}
It is a straightforward diagram check that all the isomorphisms in (\ref{TheEasyProof})
and (\ref{easyHLofDual}) are natural. In particular, that means that if we have a
morphism of complexes $P_{\bullet} \stackrel{f}{\rightarrow} Q_{\bullet}$, then
there is a commutative diagram :
$$\xymatrixcolsep{10pc}\xymatrix{
H_{-t}(Q_{\bullet}^*) \ar[r]^{H_{-t}(f^*)}_{\sim} \ar[d]^{\wr} &
H_{-t}(P_{\bullet}^*) \ar[d]^{\wr} \\
H_{t-d}(Q_{\bullet})^{\vee} \ar[r]_{H_{t-d}(f)^{\vee}}^{\sim} &
H_{t-d}(P_{\bullet})^{\vee}
}$$
\end{rem}
Finally, as an easy consequence of the above theorem, we obtain (for free!) that
$D^b_{{\mathcal A}}\left({\mathcal A}\right)$ is closed under duality $M^{\vee}=Ext^d(M, A)$.
%
%\noindent{\bf Notation.} In this section, $*$ will denote the duality in 
%${\mathcal M}FPD^{fl}(A)$, give by $M^*:=Ext^d(M, A)$. The induced  duality in 
%$D^b\left({\mathcal M}FPD^{fl}(A) \right)$ will be denoted by $\vee$.  

\bT\label{relativeDUALITY}
The category $Ch^b_{{\mathcal A}}\left({\mathcal A} \right)$ is closed under
the duality $^{\vee}$ induced by the duality $^{\vee}$ in ${\mathcal A}$. 
\eT
\pf Suppose $M_{\bullet}$ is a complex in 
$Ch^b_{{\mathcal A}}\left({\mathcal A}\right)$. 
Without loss of generality, we assume $M_{\bullet}$ is supported on 
$[n, 0]$. 
Each component $M_i$ has a projective resolution of length $d$, and putting them
together with the induced maps, we get a double complex $L_{\bullet \bullet}$, as
in the left figure below :\\
{\scalefont{0.2}
\begin{tikzpicture}[font=\tiny]
\begin{scope}
\matrix (m) [matrix of math nodes, row sep=1.5em,column sep=0.6em, text height=1.5ex, text depth=0.25ex]
{  & 0 & \vdots  &  \vdots & \cdots & \vdots & \vdots & 0  \\
L_{\bullet \bullet}= & 0 &  L_{1 n}   &  L_{1 (n-1)}  & \cdots  & L_{1 1}    & L_{1 0}   & 0  \\ 
& 0 &  L_{0 n}   &  L_{0 (n-1)}  & \cdots  & L_{0 1}  & L_{0 0}   & 0  \\ 
M_{\bullet}=& 0 &  M_n  &  M_{n-1}  & \cdots  & M_{1}  & M_0   &0\\};
\path[->]
(m-1-2) edge (m-1-3)
        edge  (m-2-2)
(m-1-3) edge  (m-1-4)
        edge  (m-2-3)
(m-1-4) edge  (m-1-5)
        edge  (m-2-4)
(m-1-5) edge  (m-1-6)
        edge  (m-2-5)
(m-1-6) edge  (m-1-7)
        edge  (m-2-6)
(m-1-7) edge  (m-1-8)
        edge  (m-2-7)
(m-2-2) edge  (m-2-3)
        edge  (m-3-2)
(m-2-3) edge  (m-2-4)
        edge  (m-3-3)
(m-2-4) edge  (m-2-5)
        edge  (m-3-4)
(m-2-5) edge  (m-2-6)
        edge  (m-3-5)
(m-2-6) edge  (m-2-7)
        edge  (m-3-6)
(m-2-7) edge  (m-2-8)
        edge  (m-3-7)
(m-3-2) edge  (m-3-3)
(m-3-3) edge  (m-3-4)
(m-3-4) edge  (m-3-5)
(m-3-5) edge  (m-3-6)
(m-3-6) edge  (m-3-7)
(m-3-7) edge  (m-3-8)
(m-4-2) edge  (m-4-3)
(m-4-3) edge  (m-4-4)
(m-4-4) edge  (m-4-5)
(m-4-5) edge  (m-4-6)
(m-4-6) edge  (m-4-7)
(m-4-7) edge  (m-4-8);
  \path[->>]
(m-3-2) edge  (m-4-2)
(m-3-3) edge  (m-4-3)
(m-3-4) edge  (m-4-4)
(m-3-5) edge  (m-4-5)
(m-3-6) edge  (m-4-6)
(m-3-7) edge  (m-4-7);
\end{scope}

\begin{scope}[xshift=7cm]
\matrix (m) [matrix of math nodes, row sep=1.5em,column sep=0.6em, text height=1.5ex, text depth=0.25ex]
{  & 0 & \vdots  &  \vdots & \cdots & \vdots & \vdots & 0  \\
L'_{\bullet \bullet}= & 0 &  L_{(d-1) 0}^*   &  L_{(d-1) 1}^*   & \cdots  & L_{(d-1) (n-1)}^*     & L_{(d-1) n}^*    & 0  \\ 
& 0 &  L_{d 0}^*    &  L_{d 1}^*   & \cdots  & L_{d (n-1)}^*   & L_{d n}^*    & 0  \\ 
M_{\bullet}^{\vee}=& 0 &  M_0^{\vee}  &  M_1^{\vee}  & \cdots  & M_{n-1}^{\vee}  & M_n^{\vee}   &0\\};
\path[->]
(m-1-2) edge (m-1-3)
        edge  (m-2-2)
(m-1-3) edge  (m-1-4)
        edge  (m-2-3)
(m-1-4) edge  (m-1-5)
        edge  (m-2-4)
(m-1-5) edge  (m-1-6)
        edge  (m-2-5)
(m-1-6) edge  (m-1-7)
        edge  (m-2-6)
(m-1-7) edge  (m-1-8)
        edge  (m-2-7)
(m-2-2) edge  (m-2-3)
        edge  (m-3-2)
(m-2-3) edge  (m-2-4)
        edge  (m-3-3)
(m-2-4) edge  (m-2-5)
        edge  (m-3-4)
(m-2-5) edge  (m-2-6)
        edge  (m-3-5)
(m-2-6) edge  (m-2-7)
        edge  (m-3-6)
(m-2-7) edge  (m-2-8)
        edge  (m-3-7)
(m-3-2) edge  (m-3-3)
(m-3-3) edge  (m-3-4)
(m-3-4) edge  (m-3-5)
(m-3-5) edge  (m-3-6)
(m-3-6) edge  (m-3-7)
(m-3-7) edge  (m-3-8)
(m-4-2) edge  (m-4-3)
(m-4-3) edge  (m-4-4)
(m-4-4) edge  (m-4-5)
(m-4-5) edge  (m-4-6)
(m-4-6) edge  (m-4-7)
(m-4-7) edge  (m-4-8);
  \path[->>]
(m-3-2) edge  (m-4-2)
(m-3-3) edge  (m-4-3)
(m-3-4) edge  (m-4-4)
(m-3-5) edge  (m-4-5)
(m-3-6) edge  (m-4-6)
(m-3-7) edge  (m-4-7);
\end{scope}
\end{tikzpicture}
}
Dualizing, $L'_{\bullet \bullet}$ gives a similar resolution of $M_{\bullet}^{\vee}$,
as shown on the right above (note that there are sign conventions on the differentials
of the complexes here, in particular $M_{\bullet}^{\vee}$ acquires a $(-1)^d$ factor on
its differentials. We refer to \cite{W} for the conventions for total complexes.) \\
Now, the total complexes give quasi-isomorphisms 
$$
Tot(L_{\bullet \bullet}) \lra M_{\bullet}, \quad 
Tot(L'_{\bullet \bullet}) \sur M_{\bullet}^{\vee}. 
$$
So, $H_i(Tot(L_{\bullet \bullet})) \in \A$
for all $i\in {\mathbb Z}$. By (\ref{easyHLofDual}), 
$H_i(Tot(L_{\bullet \bullet}^*)) \in {\mathcal A}$
for all $i\in {\mathbb Z}$. Now after translating $Tot(L_{\bullet \bullet})^*$
$d$ components to the left, we observe that it is actually chain homotopy equivalent
to $Tot(L'_{\bullet \bullet})$ and so we have $H_i(Tot(L'_{\bullet \bullet})) \in \A$.
Finally, the above quasi-isomorphism yields that 
$H_i(M_{\bullet}^{\vee}) \iso H_i(Tot(L'_{\bullet \bullet})) \in \A$. This completes the proof.
\pic $\eop$

\section{Definitions of Witt Groups}\label{WittGroups}
In this section, we define Witt groups of the categories we work with. In particular,
we extend the definition of Witt groups from triangulated categories
with duality to their additive subcategories which are closed under orthogonal sums,
translations and isomorphisms. Since it is possible that there is cause for confusion
about translation, we start by clearing the air.
\begin{defin}\label{unsignedD}{\rm
In all the categories of complexes, there are two possible translations,
$T_u$ and $T_s$. The complex $T_uP_{\bullet}$ is defined as 
$(T_u P_{\bullet})_i = P_{i-1}$ and 
$\partial(T_u P_{\bullet})_i = \partial(P_{\bullet})_{i-1}$.
The complex $T_s P_{\bullet}$ is defined as 
$(T_s P_{\bullet})_i = P_{i-1}$ and 
$\partial(T_s P_{\bullet})_i = -\partial(P_{\bullet})_{i-1}$.

{\bf {Note that $T_s$ seems to be the "standard" translation in literature and that
is always the translation we use on any category of complexes.}}

However, given a duality $*$ on such a category (e.g. $D^b_{\A}(\A)$ and $D^b_{\A}(\P(A))$),
there are shifted dualities, $T_s^n \smallcirc{0.7} *$ and $T_u^n \smallcirc{0.7} *$. 
%There are natural equivalences which allow us to go from the signed to the unsigned duality.
We work with the unsigned duality $T_u^n \smallcirc{0.7} *$ until we reach section 
\ref{shiftedWSec}.
Note however that $H_i(T_s^n P_{\bullet}^*) = 
H_i(T_u^n P_{\bullet}^*)$ and so much of what we will say is independent of the chosen duality.
% where we pass to the other one because the Witt groups $W^i$ are defined w.r.t. the signed duality. 
}
\end{defin} 
\begin{rem}\label{rmk2}
We quickly review the situation for the categories $\A$ and ${\mathcal R}$.
First note that both of these categories are exact categories with duality and so
the Witt groups are defined as in \cite{QSS}.

The functor $\iota$ induces duality preserving  equivalences
$$
\iota : \left({\mathcal A}, ^{\vee}, \tilde{\varpi}\right)
\lra 
\left({\mathcal R}(A), T^d_u \smallcirc{0.5} *, \varpi \right), \quad 
\iota : \left({\mathcal A}, ^{\vee}, -\tilde{\varpi}\right)
\lra 
\left({\mathcal R}(A), T^d_u \smallcirc{0.5} *, -\varpi \right) 
$$
of categories which then yield isomorphisms of the corresponding Witt groups
$$
W(\iota) : W\left({\mathcal A}, ^{\vee}, \tilde{\varpi}\right)
\iso
W\left({\mathcal R}(A), T^d_u \smallcirc{0.5} *, \varpi \right),
$$
$$
W(\iota) : W\left({\mathcal A}, ^{\vee}, -\tilde{\varpi}\right)
\iso 
W\left({\mathcal R}(A), T^d_u \smallcirc{0.5} *, -\varpi \right) 
$$
\end{rem}

Finally, we get to our definitions of the Witt group.
Given an exact category ${\mathcal E}$, its derived category will be denoted by 
$D^b({\mathcal E})$. For a subcategory ${\mathcal C}$, 
$D^b_{{\mathcal C}}({\mathcal E})$ will denote the full subcategory of 
$D^b({\mathcal E})$ consisting of complexes with homologies in ${\mathcal C}$.

The derived categories which will be used in this article include :
$$
D^b_{{\mathcal A}}({\mathcal P}(A)) \subseteq D^b_{fl}({\mathcal P}(A))  \subseteq D^b({\mathcal P}(A)) ~~~~~~~~~~ D^b_{{\mathcal A}}({\mathcal A}) \subseteq D^b(\A).
$$

We have the following diagram of subcategories and functors: 
\begin{equation}\label{diag1}
\diagram 
{\mathcal A} \ar[d]_{\iota}\ar[r]^{\mu\qquad} \ar[dr]^{\zeta}  
& D^b_{{\mathcal A}}\left({\mathcal A} \right) 
\ar@{^(->}[r]^{\quad\nu}\ar[d]^{\alpha}
& D^b\left({\mathcal A} \right)\ar[d]^{\beta} & \\
{\mathcal R}(A) \ar@{^(->}[r]^{\mu'}
& D^b_{{\mathcal A}}\left({\mathcal P}(A)\right)
\ar@{^(->}[r]^{\nu'} & D^b_{fl}\left({\mathcal P}(A)\right).
\enddiagram 
\end{equation}
Here $\mu(M)$ is the complex concentrated at degree zero. The functor $\iota(M)$ is
obtained by making a choice of projective resolution of length $d$ and then defining
$M^{\vee} = H_0((\zeta(M))^*)$. The functors $\alpha$ and $\beta$ are essentially
induced by these ones, by taking the total complex (take a look at the proof of
(\ref{relativeDUALITY})).

We now move on towards the definitions of the Witt groups of the categories
$D^b_{{\mathcal A}}({\mathcal A})$ and $D^b_{{\mathcal A}}({\mathcal P}(A))$.
We once again remind the reader that this definition relies on the definitions
in \cite{TWGI}.

%%%%%%%%%%%%%%%%%%%%%%%%%%%%We might want to include those definitions.

\begin{defin} \label{defWittGr}
{\rm Let $\delta = \pm 1$.
 Suppose $K:=(K,\#, \delta, \varpi)$ is a triangulated category with
translation $T$ and  $\delta-$duality $\#$. 
Suppose $K_0$ is a full subcategory of $K$ that is closed under isomorphism, translation
and orthogonal sum. We abuse notation and denote $K_0 :=(K_0,\#, \delta,\varpi)$ in order to
keep track of the duality and canonical isomorphism in use.
  
\bE
\item Define the Witt monoid of $MW(K_0)$ to be the submonoid 
$$
MW(K_0)=\{(P,\varphi)\in MW(K): P \in Ob(K_0) \}.
$$

\item \label{aboutNeutral} A symmetric space
$(P, \varphi) \in MW(K_0)$ will be called a {\bf neutral space} in $MW(K_0)$
if it has a Lagrangian $(L, \alpha, w)$ in $MW(K)$ such that $L, L^{\#} \in Ob(K_0)$.

\item  
Let $NW(K_0)$ be the submonoid of $MW(K_0)$ generated by
the isometry classes of neutral spaces in $K_0$.  

\item Define the Witt group 
$$
W(K_0):=\frac{MW(K_0)}{NW(K_0)}.
$$ 
Note $(Q,\chi)\in MW(K_0) \Lra (Q,-\chi)\in MW(K_0)$. It is easy to check that
$(Q,\chi)\perp (Q,-\chi)\in NW(K_0)$. So, $W(K_0)$ has a group structure. We use
this definition in the context of derived categories of exact categories with duality. 
%
%\item The Witt group of of the derived category 
%$D^b\left({\mathcal M}FPD^{fl}(A)\right)$ (which is already define in the literature), will be denoted by 
%$W\left(D^b\left({\mathcal M}FPD^{fl}(A)\right)\right)$.
%\item With $K_0=D^b_{{\mathcal M}FPD^{fl}}\left({\mathcal P}(A)\right)$ 
%and $\#=T^d*$ we define the Witt group  
%$$ 
%W\left(D^b_{{\mathcal M}FPD^{fl}}\left({\mathcal P}(A)\right), \#\right). 
%$$ 
\item
Let $(\C, ^{\vee}, \varpi)$ be an exact subcategory with duality in an ambient abelian
category $\C'$ and let $\D$ be any subcategory of $\C'$ closed under orthogonal sum.
Let $K_0=D^b_{{\D}}({\C})$~(\ref{notaCate}). Then with the induced duality and natural
isomorphism, the Witt group $W\left(D^b_{\D}(\C), ^{\vee}, \delta, \varpi \right)$ is defined
as above.
\item Accordingly, with $T= T_u, T_s$, the Witt groups 
$$
W\left(D^b_{{\mathcal A}}\left({\mathcal P}(A)\right), T^n \smallcirc{0.5} *, \pm 1, \pm \varpi  \right), ~~~
W\left(D^b_{\A}\left( \A \right), T^n \smallcirc{0.5} *, \pm 1, \pm \tilde{\varpi} \right) 
$$
are defined. 
%\item {\bf Notations:} We denote 
%$$
%D^b_{{\mathcal A}}\left({\mathcal P}(A)\right)^{\pm}_u:=
%\left(D^b_{{\mathcal A}}\left({\mathcal P}(A)\right), T_u, *,  1,\pm \varpi \right), ~~
%D^b_{{\mathcal A}}\left({\mathcal P}(A)\right)^{\pm}_s:=
%\left(D^b_{{\mathcal A}}\left({\mathcal P}(A)\right), T_s, *,  1,\pm \varpi  \right). 
%$$
%Similalarly, define
%$$
%\left(D^b_{{\mathcal A}}\left(\P(A) \right)^{\pm}_u\right):=\left(D^b_{{\mathcal A}}\left({\mathcal A}\right), T_u, *,  1,\pm \varpi \right), ~~
%D^b_{{\mathcal A}}\left({\mathcal A}\right)^{\pm}_s
%:=\left(D^b_{{\mathcal A}}\left({\mathcal A}
%\right), T_s, *,  1,\pm \varpi  \right). 
%$$
\eE
}\end{defin} 
\section{Isomorphisms of Witt Groups} \label{WDIAGRAM} 
All the functors above induce homomorphisms of Witt groups.
As always, {\bf $A$  denotes a Cohen-Macaulay ring with $\dim A_m=d \geq 2$,
for all maximal ideals $m$ of $A$}.
Let $D^b_{{\mathcal A}}
\left({\mathcal A}, ^{\vee}, \pm\tilde{\varpi} \right)$, 
 $D^b \left({\mathcal A}, ^{\vee}, \pm\tilde{\varpi}  \right)$
denote the duality structure, respectively, on 
 $D^b_{{\mathcal A}} \left({\mathcal A}\right)$ 
 and $D^b\left({\mathcal A}\right)$ induced by 
$\left({\mathcal A}, ^{\vee}, \pm\tilde{\varpi}\right)$ and
$$
D^b_{{\mathcal A}}\left( \P(A) \right)^{\pm}_u:=
\left(D^b_{{\mathcal A}}({\mathcal P}(A)), T^d_u \smallcirc{0.5} *, 1, \pm \varpi \right).
$$ 
%%%%%%%%%%%%%%%%%%%%%%%%%We probably need to say why the maps in the diagram below exist.
Recall the functors in the diagram \eqref{diag1}, it is clear that the functors $\mu,\nu,\mu'$
and $\gamma$ preserve dualities. For $\iota$, this is left to the reader as a diagram chase
(but note the definition of $M^{\vee}$ after that same diagram). Essentially
the same proof also gives us that $\zeta, \alpha$ and $\beta$ are duality preserving. This
being done, we can talk about the corresponding maps of the Witt groups. The goal of this
section is to establish the following diagram of homomorphisms of Witt groups :
\begin{equation}\label{diag2}
\xymatrix{
W\left({\mathcal A}, ^{\vee}, \pm\tilde{\varpi}\right) 
\ar[rd]^{\sim}_(0.4){W(\zeta)} 
\ar[d]_{W(\iota)}^{\wr} \ar[r]_{\sim}^{W(\mu)}
& W\left(D^b_{{\mathcal A}}
\left({\mathcal A}, ^{\vee}, \pm\tilde{\varpi} \right) \right) 
\ar[r]_{\sim}^{W(\nu)} \ar[d]^{\wr}_{W(\alpha)} 
&W\left(D^b\left({\mathcal A}, ^{\vee}, \pm\tilde{\varpi}\right) \right) \\
W\left({\mathcal R}(A), T^d_u*, \pm\varpi \right) \ar[r]_{\sim}^{W(\gamma)}
& W\left(D^b_{{\mathcal A}}\left(\P(A) \right)^{\pm}_u\right)
%D^b_{{\mathcal A}}\left({\mathcal P}(A), T_u, T^d_u*, 1, 
%\pm \varpi\right)\right)
&\\ %\ar[r] & W^d\left(D^b_{fl}\left({\mathcal P}(A)\right)\right)
}
\end{equation}
Note that we already know that $W(\iota)$ is an isomorphism (\ref{rmk2}) and further,
by \cite[Theorem 4.3]{TWGII}, $W(\nu \circ \mu)$ is an isomorphism. The proof that
$$
W(\mu): W\left({\mathcal A}, ^{\vee}, \pm \tilde{\varpi} \right) \lra
W\left(D^b_{{\mathcal A}}
\left({\mathcal A}, ^{\vee}, \pm \tilde{\varpi} \right)\right) 
%W\left(D^b_{{\mathcal M}FPD^{fl}}
%\left({\mathcal M}FPD^{fl}(A) \right)  \right) 
$$
are  isomorphisms follows from more abstract "general nonsense" which we prove in
(\ref{midWIsoF}) as part of the appendix (\ref{appA}).
Since $W(\nu \circ \mu)$ and $W(\mu)$ are isomorphisms, it is clear that so is $W(\nu)$.
The main result of this section is that $W(\zeta)$ is an isomorphism. That being established,
it is clear that $W(\gamma)$ and $W(\alpha)$ are isomorphisms.

For the rest of this section, we use the notation $\#:=T^d_u \smallcirc{0.5} *$. 
First, we establish the following regarding the structure of
symmetric forms. 
%$D^b_{{\mathcal M}FPD^{fl}}\left({\mathcal P}(A)\right)$. 
\begin{lemma}\label{stOfForms}
Suppose $\eta:X_{\bullet} \iso X^{\#}_{\bullet}$ is a symmetric form in 
$\left(D^b_{{\mathcal A}}\left(\P(A) \right)^{\pm}_u\right)$, such that   
%$D^b_{{\mathcal M}FPD^{fl}}\left({\mathcal P}(A)\right)$
%(with duality $\#=T^d*$). 
$$
H_{-m}(X_{\bullet}) \neq 0, \quad and \quad H_{i}(X_{\bullet})=0 \quad for
~all\quad i<-m. 
$$
Then, there is a complex $P_{\bullet}$ in $Ch^b_{{\mathcal A}}\left({\mathcal 
P}(A)\right)$ 
%of projective modules 
and 
a quasi-isomorphism $\varphi:P_{\bullet} \rightarrow P_{\bullet}^{\#}$ 
such that 
%%%%%%%%%%%%%%%%%%%%%%%%%%%%%%%%%%%%%%%%%%%%Do you wanna replace quasi by homotopy equiv.?
\bE
\item $(P_{\bullet}, \varphi)$ is isometric to $(X_{\bullet},\eta)$
in 
$\left(D^b_{{\mathcal A}}\left(\P(A) \right)^{\pm}_u\right)$.  
%$D^b_{{\mathcal M}FPD^{fl}}\left({\mathcal P}(A)\right)$.
\item $P_{\bullet}$ is supported on $[m+d,-m]$.
\item $H_{-m}(P_{\bullet}) \neq 0$. 
\eE 
\end{lemma}
\pf Recall from (\ref{GoodAtLeft}) that since $H_{-m}(X_{\bullet}) \neq 0$, 
$X_{\bullet}$ has length at least $d$. By duality, we conclude that $m \geq 0$.
By definition there is a complex $P_{\bullet}$ of projective modules
and a quasi-isomorphism  $t:P_{\bullet}\lra X_{\bullet}$, a chain complex morphism
$\varphi_0: P_{\bullet}\lra X^{\#}_{\bullet}$ such that $\eta=\varphi_0 t^{-1}$.
Then, $\varphi=t^{\#}\varphi_0= t^{\#}\eta t$ is
a symmetric form on
$P_{\bullet}$, and $(X_{\bullet}, \eta)$ is isometric to $(P_{\bullet}, \varphi)$. 
By including enough zeros on the two tails, we can assume $P_{\bullet}$ 
is supported on $[n+d,-n]$, for some $n \geq m$.
If $m=n$ there is nothing to prove. So, assume $n>m$. 
We have, $H_{-n}(P_{\bullet})= 0$. 
Inductively, we will cut down the support to $[m+d, m]$.
We write $ \varphi : P_{\bullet}\lra P_{\bullet}^{\#}$ as follows
{\scalefont{0.7}
$$
\diagram
\ar[r] & 0 \ar[r] \ar[d]& 
P_{n+d} \ar[d]_{\varphi_n}\ar[r]^{\partial_{n+d}}
&P_{(n-1)+d} \ar[d]_{\varphi_{(n-1)+d}} \ar[r]^{\partial_{(n-1)+d}} & \cdots \ar[r]^{\partial_{-(n-2)}}
& P_{-(n-1)}\ar[r]^{\partial_{-(n-1)}} \ar[d]_{\varphi_{-(n-1)}} 
& P_{-n} \ar[d]^{\varphi_{-n}}\ar[r] & 0 \ar[d] \ar[r] & \\
\ar[r] & 0 \ar[r] & P_{-n}^* \ar[r]_{\partial_{-(n-1)}^*}&P^*_{-(n-1)}\ar[r]^{\partial_{-(n-2)}^*}
&\cdots \ar[r]_{\partial_{(n-1)+d}^*} & P_{(n-1)+d}^*\ar[r]_{\partial_{n+d}^*} & P_{n+d}^* \ar[r] & 0  \ar[r] & \\ 
\enddiagram 
$$
}
where $P_i$ are finitely generated projective $A-$modules. Since $n>m, 
H_{-n}(P_{\bullet}) \cong H_{-n}(P_{\bullet}^*) \cong 0$.
So, $\partial_{-(n-1)}$ and $\partial_{n+d}^*$ are both split surjections.
Thus there are homomorphisms $\epsilon_{-n}: P_{-n}\lra P_{-(n-1)}$ and
$\epsilon_{n+d}^*: P_{n+d}* \lra P_{n-1+d}$ such that 
$\partial_{-(n-1)} \circ \epsilon_{-n} = Id$ and
$\partial_{n+d}^* \circ \epsilon_{n+d}^* = Id$.
Hence, $Z_{-(n-1)}$ and $\frac{P_{n-1+d}}{P_{n+d}}$ are projective modules.
Note that since $d \geq 2$, by (\ref{GoodAtLeft}), $\frac{P_{n-1+d}}{P_{n+d}} = B_{n-2+d}$.
Further, we obtain splittings $\sigma_{-(n-1)} : P_{-(n-1)} \sur Z_{-(n-1)}$ and
$\sigma_{(n-2)+d} : B_{(n-2)+d} \hookrightarrow P_{(n-1)+d}$. This gives us a shorter
complex $Q_{\bullet}$, naturally chain homotopic to $P_{\bullet}$ and an induced symmetric
form on $Q_{\bullet}$ :
{\scalefont{0.8}
$$
\xymatrix@C-=0.5cm{
Q_{\bullet} : ~~~~ \ar[r] \ar[d] & 0 \ar[r] \ar[d]
& 0 \ar[d]\ar[r] & B_{(n-2)+d} \ar[d]^{\sigma_{(n-2)+d}} \ar@{^(->}[r] & P_{(n-2)+d}~ \cdots \ar[r]^{\partial} \ar@{=}[d] & Z_{-(n-1)}\ar[r] \ar@{^(->}[d] & 0 \ar[d] \ar[r] & 0 \ar[d] \ar[r] & \\
P_{\bullet} : ~~~~ \ar[r] \ar[d]_{\varphi} & 0 \ar[r] \ar[d] & P_{n+d} \ar[d]_{\varphi_{n+d}} 
\ar[r]^{\partial_{n+d}} & P_{(n-1)+d} \ar[d]_{\varphi_{(n-1)+d}} 
\ar[r]^{\partial} & P_{(n-2)+d}~ \cdots \ar[r]^{\partial} \ar[d]_{\varphi_{(n-2)+d}}
 & P_{-(n-1)} \ar[r]^{\partial} \ar[d]_{\varphi_{-(n-1)}}  & P_{-n} \ar[d]^{\varphi_{-n}} 
\ar[r] & 0 \ar[d] \ar[r] & \\
P_{\bullet}^{\#} : ~~~~ \ar[r] \ar[d] & 0 \ar[r] 
 \ar[d] & P_{-n}^* \ar[d] \ar[r]_{\partial_{-(n-1)}^*} &
P^*_{-(n-1)} \ar@{->>}[d] \ar[r]_{\partial_{-(n-2)}^*} & P^*_{-(n-2)} \cdots \ar@{=}[d] \ar[r]_{\partial^*} & P_{(n-1)+d}^* \ar[d]^{\sigma_{(n-2)+d}^*} \ar[r]_(.7){\partial_{n+d}^*} 
& P_{n+d}^* \ar[d] \ar[r] & 0 \ar[d] \ar[r] & \\ 
Q_{\bullet}^{\#} : ~~~~ \ar[r] & 0 \ar[r] & 0 \ar[r] & Z^*_{-(n-1)} \ar[r] 
& P^*_{-(n-2)} \cdots \ar[r]_{\partial^*} & B_{(n-2)+d}^* \ar[r] & 0 \ar[r] & 0 \ar[r] &
}
$$
}
Calling this map $\varphi'$, $(Q_{\bullet}, \varphi')$ is obviously isometric to
$(P_{\bullet}, \varphi)$ and hence to the original form $(X_{\bullet}, \eta)$.
Since $Q_{\bullet}$ is supported on $[(n-1)+d, -(n-1)]$ induction finishes the proof.
$\eop$ 

Since $\zeta$ is given by composing $\nu$ and $\alpha$ and there are maps $W(\nu)$ and
$W(\alpha)$, it is clear that $W(\zeta)$ is well-defined. However, we give an explicit
proof which might also be somewhat illuminating considering the unsaid details about
why duality-preserving functors induce maps of Witt groups. The proof essentially follows
the proof in \cite[2.11]{DWG}.
\bT\label{WzetaWellD}
The functor $\zeta$ induces a well defined homomorphism
$$
W(\zeta): 
W\left({\mathcal A}, ^{\vee}, \pm\tilde{\varpi} \right)
\lra W\left(\left(D^b_{{\mathcal A}}\left(\P(A) \right)^{\pm}_u\right) \right). 
$$
\eT 
\pf We will only prove 
$$
W(\zeta): 
W \left({\mathcal A}, ^{\vee}, \tilde{\varpi} \right)
\lra W\left(D^b_{{\mathcal A}}\left({\mathcal A}\right)^{+}_u \right). 
$$
is well defined and 
the case of skew dualities follows similarly. 
%%%%%%%%%%%%%%%%%%%%%%%%%%%%%%%%%%%%%%%%%%%%%%%%%%%%%%%%%%%%%%%%%%%%%%%%%
It is clear that $\zeta$ defines a well-defined map from
$MW\left({\mathcal A}, ^{\vee}, \tilde{\varpi} \right)$ to 
$MW\left(D^b_{{\mathcal A}}\left({\mathcal A}\right)^{+}_u \right)$ since
projective maps of modules can be lifted to a chain complex map of their resolutions
(note that though the lift is not unique, it is unique upto homotopy and so gives the
same morphism in $D^b_{\A}(\P(A))$). So we need to check that the image of a neutral space
in $MW\left({\mathcal A}, ^{\vee}, \tilde{\varpi} \right)$  is neutral in 
$MW\left(D^b_{{\mathcal A}}\left({\mathcal A}\right)^{+}_u \right)$.

Suppose $(M, \varphi_0)$ is a neutral space in 
$\left({\mathcal A}, ^{\vee}, \tilde{\varpi} \right)$.
Let $\alpha_0: N \lra M$ be a lagrangian of $(M, \varphi_0)$.
Then
$$
\diagram 
0 \ar[r] & N \ar[r]^{\alpha_0}
& M \ar[r]^{\alpha_0^{\vee}\varphi_0} & N^{\vee} \ar[r] & 0   
\enddiagram
\quad is ~exact. 
$$
Suppose $L_{\bullet}, P_{\bullet}$ are the chosen projective resolutions of $N$ and $M$
and $\alpha : L_{\bullet} \rightarrow P_{\bullet}$ is the morphism induced from $\alpha_0$.
The above short exact sequence implies the composition
$L_{\bullet} \stackrel{\alpha}{\rightarrow} P_{\bullet} \stackrel{\alpha^{\#}\varphi}{\rightarrow}  L^{\#}_{\bullet}$ is chain homotopic to $0$ (hence the $0$ map in $D^b_{\A}(\P(A))$). Completing $\alpha$ to an exact triangle, we get a morphism of exact triangles
$$\xymatrix{
L_{\bullet} \ar[r]^{\alpha} \ar[d] & P_{\bullet} \ar[r]^{j} \ar[d]^{\alpha^{\#}\varphi} & C_{\bullet} \ar[r]^{k} \ar@{-->}[d]^{s} & T(L_{\bullet}) \ar[d] \\
0 \ar[r] & L_{\bullet}^{\#} \ar@{=}[r] & L_{\bullet}^{\#} \ar[r] & 0\\
}$$
Note that $H_0(C_{\bullet})\cong N^{\vee}$ and $\forall ~~i\neq 0\quad H_i(C_{\bullet})=0$
and so $C_{\bullet}$ is an object in $D^b_{\A}(\P(A))$. The map $s$ is actually quite easy
to describe, namely $s=(0, \alpha^{\#}\varphi): L_{n-1} \oplus P_n \lra L^*_n$ and it follows
from the above morphism of triangles (or by direct checking) that $s$ is a quasi-isomorphism.
Hence, 
$$\xymatrix{
L_{\bullet} \ar[r]^{\alpha} & P_{\bullet} \ar[r]^{\alpha^{\#}\varphi} & L_{\bullet}^{\#}
 \ar[r]^{k \circ s^{-1}} & T(L_{\bullet})\\
}$$
is an exact triangle.
Setting $w=-T^{-1}(k \circ s^{-1})$, we get an exact triangle 
$$\xymatrix{
T^{-1}(L_{\bullet}^{\#}) \ar[r]^w & L_{\bullet} \ar[r]^{\alpha} & P_{\bullet} \ar[r]^{\alpha^{\#}\varphi} & L_{\bullet}^{\#} \\
}$$
Now all we require is that $T^{-1}w^{\#} = w$.
\begin{align*}
 T^{-1}w^{\#} = w \Leftrightarrow T^{-1}w^{\#} = -T^{-1}(k \circ s^{-1})
& \Leftrightarrow (T^{-1}(k \circ s^{-1}))^{\#} = k \circ s^{-1} \\
\Leftrightarrow T({s^{-1}}^{\#} \circ k^{\#}) = k \circ s^{-1}
& \Leftrightarrow T(k^{\#}) \circ s = T(s^{\#}) \circ k.
\end{align*}
A quick physical check of the maps in question yields that the first map is
$$\xymatrixcolsep{10pc}\xymatrix{
L_{n-1} \oplus P_n \ar[r]^{ \left( \begin{array}{cc} -1 & 0 \end{array} \right)} &
L_{n-1} \ar[r]^{ \left( \begin{array}{c} 0 \\ \varphi_{n-1}^* \circ \alpha_{d-n+1} \end{array} \right)} & L_{d-n}^* \oplus P_{d-n+1}^*}$$
while the second one is
$$\xymatrixcolsep{10pc}\xymatrix{
L_{n-1} \oplus P_n \ar[r]^{ \left( \begin{array}{cc} 0 & \alpha_{d-n}^* \circ \varphi_n \end{array} \right)} &
L_{d-n}^* \ar[r]^{ \left( \begin{array}{c} -1 \\ 0 \end{array} \right)} & L_{d-n}^* \oplus P_{d-n+1}^*}$$
The matrices we thus obtain are
$$\left( \begin{array}{cc}
 0 & -\varphi_{n-1}^* \circ \alpha_{d-n+1}\\
 0 & 0 
\end{array}\right) ~~~~~~~~~~~~~~~~
\left(\begin{array}{cc}
 0 & 0  \\
\alpha_{d-n}^* \circ \varphi_n & 0 
\end{array}\right)$$
which are homotopy equivalent using 
$$
\tau=
\left(\begin{array}{cc}0 &0\\ 0& (-1)^n\varphi \end{array}\right).$$
Therefore, $(L_{\bullet}, \alpha, w)$ is a lagrangian.
Hence $W(\zeta)$ is a well defined homomorphism of groups. 
$\eop$ 

%%%%%%%%%%%%%%%%%%%%%%%%%%%%%%%%%%%%%%%%%%%%%%%%%%%%%%%%%%%%%%%%%%%%%%%%%
\vspace{5mm}
We now proceed towards (\ref{reduceL}) which proves that $W(\zeta)$ is surjective.
The main tool here is to construct a special sublagrangian and then use 
Balmer's sublagrangian construction \cite[Section 4 and Theorem 4.20]{TWGI}
to reduce the length of $(P,\varphi)$.
%%%%%%%%%%%%%%%%%%%%%%%%%%%%%%%%%%%%%%%%%%%%%%%%%%%%%%%%%%%%%%%%%%%%%%%%%
\begin{rem}\label{rmk3}
Note that using (\ref{stOfForms}), any symmetric form $(X_{\bullet}, \phi)$
in $\left(D^b_{{\mathcal A}}\left({\mathcal A}\right)^{+}_u \right)$ with $X_{\bullet}$
not acyclic can be represented by
$$
\diagram
P_{\bullet}=\ar[d]_{\varphi} \cdots~0 \ar[r] 
& P_{n+d} \ar[d]_{\varphi_n}\ar[r]^{\partial}
&P_{(n-1)+d}\ar[d]_{\varphi_{(n-1)+d}}  \ar[r]^{\partial}\ar[r] & \cdots \ar[r]^{\partial} 
& P_{-(n-1)}\ar[r]^{\partial} \ar[d]_{\varphi_{-(n-1)}} 
& P_{-n} \ar[d]^{\varphi_{-n}}\ar[r] & 0 \\
P_{\bullet}^{\#}=\cdots~0 \ar[r] & P_{-n}^* \ar[r]_{\partial^*}&P^*_{-(n-1)}\ar[r] 
&\cdots \ar[r]_{\partial^*}& P_{(n-1)+d}^*\ar[r]_{\partial^*} & P_{n+d}^* \ar[r] & 0  \\ 
\enddiagram 
$$
with $H_{-n}(P_{\bullet}) \neq 0$.
\end{rem}
\begin{lemma}\label{HnNotZero} 
Let $(P_{\bullet}, \varphi)$ be as above.
Then 
\bE
\item $H_r(P_{\bullet}) =0$ for $r=n+1, n+2, \ldots, n+d$.
\item $H_n(P_{\bullet}) \neq 0$.
\eE 
\end{lemma}
\pf The first point follows from 
(\ref{GoodAtLeft}). To prove (2), assume $H_n(P_{\bullet})=0$. 
Then, with $B_{n-1}=image(\partial_{n})$ we have an exact sequence
$$
\diagram 
0\ar[r] & P_{n+d} \ar[r] & P_{(n-1)+d} \ar[r] & \cdots \ar[r] 
& P_n \ar[r] & P_{n-1} \ar[r] & \frac{P_{n-1}}{B_{n-1}}\ar[r] & 0 
\enddiagram 
$$
Since this is a projective resolution of the last term, if follows
$$
H_{-n}(P_{\bullet}^*) = Ext^{d+1}\left( \frac{P_{n-1}}{B_{n-1}}, A \right) = 0.
$$
This is
a contradiction to $H_{-n}(P_{\bullet}) \neq 0$. \pic $\eop$

\vspace{5mm}  
Much of what follows is dependent on \cite[\S 4]{TWGI} and the interested
reader is highly encouraged to take a look at it. We recall the definition of
a sublagrangian of $(P_{\bullet}, \varphi)$ :
\begin{defin}
A sublagrangian of a symmetric form $(P_{\bullet}, \varphi)$ is a pair 
$(L_{\bullet},\alpha)$ with $L_{\bullet} \in Ob(D^b_{\A}(\P(A)))$ and
$\alpha : L_{\bullet} \rightarrow P_{\bullet}$, which satisfies that
$\alpha^{\#} \smallcirc{0.7} \varphi \smallcirc{0.7} \alpha = 0$ in $D^b_{\A}(\P(A))$.
\end{defin}
For $(P_{\bullet}, \varphi)$ as above, (\ref{HnNotZero}) tells us that 
$H_{n}(P_{\bullet}) \neq 0$ and we already know it is in $\A$. So it has a minimal
projective resolution of length $d$. Let $L_{\bullet}$ be a projective resolution
of $H_n(P_{\bullet})$ of length $d$, shifted by $n$ places, as in the diagram below.
Since $H_i(P_{\bullet}) ~\forall~i > n$ by (\ref{GoodAtLeft}), the bottom line is
a projective resolution of $\frac{P_n}{B_n}$ and so the inclusion $H_n(P_{\bullet}) \hookrightarrow \frac{P_n}{B_n}$ induces a map of complexes
$$
\diagram 
 0\ar[r] & L_{n+d} \ar[r]\ar[d]_{\nu_{n+d}}  
& L_{(n-1)+d} \ar[d]_{\nu_{(n-1)+d}} \ar[r] & \cdots \ar[r] 
& L_n \ar[d]_{\nu_{n}} \ar[r] \ar[r] & 0\\
0\ar[r] & P_{n+d} \ar[r] & P_{(n-1)+d} \ar[r] & \cdots \ar[r] 
& P_n \ar[r] & 0.\\
\enddiagram 
$$
Note that since the composition $H_n(P_{\bullet}) \hookrightarrow \frac{P_n}{B_n} 
\twoheadrightarrow B_{n-1}$ is $0$, we get a chain complex morphism 
$\nu: L_{\bullet} \lra P_{\bullet}$.
\begin{lemma}\label{leftSubLag}
With the notations as above, for $n > 0, (L_{\bullet}, \nu)$ defines a sublagrangian
of $(P_{\bullet}, \varphi)$.
\end{lemma} 
\pf 
Let $\alpha= \nu^{\#}\varphi \nu$ is as follows ({\it the first line indicates
the degrees}):
{\scalefont{0.7}
$$
\diagram 
    &   &  {n+1}\ar@{-->}[d] 
  & n \ar@{-->}[d]  & n-1 \ar@{-->}[d] &  & -n\ar@{-->}[d] & &\\
\ar[r]  & L_{n+2} \ar[d]_{\alpha_{n+2}} \ar[r]^{\partial} 
 &  L_{n+1}\ar[d]_{\alpha_{n+1}} \ar[r]^{\partial} 
 & L_n \ar[d]_{\alpha_{n}} \ar[r] & 0\ar[d] \ar[r] & \cdots & 0 \ar[r] \ar[d] & 0 \ar[r] &\\
 \ar[r] & L_{d-(n+2)}^* \ar[r]_{\partial_{d-n-1}^*} &  L_{d-(n+1)}^* \ar[r]_{\partial_{d-n}^*} 
 & L_{d-n}^* \ar[r]_{\partial_{d-n+1}^*} & L_{d-n+1}^* \ar[r]_{\partial^*} & \cdots \ar[r] 
 & L_{n+d}^* \ar[r] & 0 \ar[r] &\\
 \enddiagram 
 $$}
$L^{\#}$ is exact at all degrees except $-n$. 
Since $n > 0, H_i(L^{\#}) = 0 \forall i \geq n$. Hence, 
$image(\alpha_n) \subseteq  ker({\partial}_{d-n+1}^*) = image({\partial}_{d-n}^*)$.
So, $\alpha_n$ lifts to a homomorphism $h_n:L_n\lra L_{d-(n+1)}^*$, i.e. ${\partial_{d-n}}^*h_n=\alpha_n$. So  ${\partial_{d-n}}^*(\alpha_{n+1}-h_n \partial_{n+1})=0$.
Now we can inductively define a homotopy $h_r: L_r \lra L^*_{d-(r+1)}$ so that 
$\alpha$ is homotopic to zero. \pic $\eop$

\vspace{5mm} 
We intend to apply the sulagrangian construction of Balmer 
\cite[Theorem 4.20]{TWGI} to $\nu$. Since $D^b_{\A}(\P(A))$ is not a triangulated
category (in particular not closed under cones), we need to reprove some of the
results in \cite[Theorem 4.20]{TWGI}. The main (and only) thing we have to keep
track of is that in all the constructions, our objects remain within the category
$D^b_{\A}(\P(A))$. We start by checking that the cone of $\nu$ constructed in
(\ref{leftSubLag}) is an object of $D^b_{\A}(\P(A))$.
\begin{lemma}\label{coneOfsubL}
With the notations of (\ref{leftSubLag}),
let $N_{\bullet}$ be the cone of $\nu$. Then,
\bE
\item 
$N_{\bullet}$ is in 
$D^b_{{\mathcal A}}\left({\mathcal P}(A) \right)$. 
\item The homologies are given by
$$
H_i(N_{\bullet})=\left\{ 
\begin{array}{ll}
H_i(P_{\bullet}) & if~ n> i\geq -n\\
0&otherwise\\ \end{array}  
\right.
$$
\item $N_{\bullet}$ is supported on $[n+d+1, -n]$.  
\eE
\end{lemma} 
\pf The last point is obvious from the construction of the cone and (1) follows from (2).
We prove (2). We have the exact triangle
$$
\diagram 
T^{-1}N_{\bullet} \ar[r] & L_{\bullet} \ar[r]^{\nu} 
& P_{\bullet} \ar[r] & N_{\bullet}.  
\enddiagram 
$$
By construction, $H_{n}(L_{\bullet})\iso H_{n}(P_{\bullet})$ 
and $H_{i}(L_{\bullet})=0$ for all $i\neq n$. 
The long exact sequence of homologies 
{\scalefont{0.7}
$$\xymatrixcolsep{1.3pc}\xymatrix{
\cdots H_{n+2}(N_{\bullet})\ar[r] & 0\ar[r] & 0\ar[r] & H_{n+1}(N_{\bullet})\ar[r] & H_{n}(L_{\bullet})\ar[r]^{\sim} & H_{n}(P_{\bullet})\ar[r] &  \\
H_{n}(N_{\bullet})\ar[r] & 0 \ar[r] & H_{n-1}(P_{\bullet})\ar[r] & H_{n-1}(N_{\bullet})\ar[r]
& 0 \ar[r] & \cdots \cdots \ar[r] & \\ 
0 \ar[r] & H_{-n}(P_{\bullet}) \ar[r] & H_{-n}(N_{\bullet}) \ar[r] & 0 \ar[r] & 0 \ar[r]
& H_{-(n+1)}(N_{\bullet}) \ar[r] & 0 \cdots
}
$$}
establishes (2) and hence the lemma.
$\eop$
%that $H_i(N_{\bullet})=0$ for $i\geq n,~ i<-n$
%and otherwise $H_i(N_{\bullet}) \cong H_{i}(P_{\bullet})$.
%So, $N_{\bullet} \in 
%D^b_{{\mathcal M}FPD^{fl}}\left({\mathcal P}(A) \right)$. 
%
%\pic $\eop$ 

\vspace{5mm} 

Now we consider the dual $N_{\bullet}^{\#}$ of the cone of $\nu$.
\begin{lemma}\label{homoOfM} 
With the same notations as above (in (\ref{leftSubLag})),
consider the following morphism of exact triangles :
$$\xymatrixcolsep{5pc}\xymatrix{
T^{-1} N_{\bullet}\ar[r]^{\nu_0}\ar[d]_{T^{-1}\mu_0}
& L_{\bullet} \ar[r]^{\nu}\ar[d]_{\mu_0} 
& P_{\bullet} \ar[r]^{\nu_2}\ar[d]_{\varphi} 
& N_{\bullet}\ar[d]_{\mu_0^{\#}}\\   
T^{-1} L_{\bullet}^{\#} \ar[r]_{T^{-1}\nu_0^{\#}} 
& N_{\bullet}^{\#} \ar[r]_{\nu_2^{\#}} 
& P_{\bullet}^{\#}  \ar[r]_{\nu^{\#}} & L_{\bullet}^{\#} \\ 
}
$$ 
(refer \cite[4.3]{TWGI}...the existence of $\mu_0$ is assured by combining axioms (TR1) and (TR3) of triangulated categories and using that $2$ is invertible.) \\
Then,
\bE
\item $N_{\bullet}^{\#}$ is in 
$D^b_{{\mathcal A}}\left({\mathcal P}(A) \right)$. 
\item $N_{\bullet}^{\#}$ is supported on $[n+d, -(n+1)]$.

\item \label{mu0iso} $\mu_0$ induces an isomorphism
 of the $n^{th}-$homology 
 $$
 H(\mu_0): H_n(L_{\bullet}) \iso H_n(N_{\bullet}^{\#}).
 $$
\item \label{otherH0}
$$
H_i(N_{\bullet}^{\#}) \cong \left\{ 
\begin{array}{ll}
H_i(P_{\bullet}^{\#}) & if~ n\geq  i> -n\\
0&otherwise\\ \end{array}  
\right.
$$
\eE 
\end{lemma} 
\pf  
(1) follows directly from (\ref{TheEasyProof}). (2) follows because by (\ref{coneOfsubL})
$N_{\bullet}$ is supported on $[(n+1)+d,-n]$. For (\ref{mu0iso}), notice that the only
nonzero homology of $L_{\bullet}^{\#}$ is at degree $-n$. Since $n > 0$, the long exact
homology sequence of the second triangle gives us 
$$H(\nu_2^{\#}): H_n(N_{\bullet}^{\#}) \iso H_n(P_{\bullet}^{\#}).$$
By choice of $\nu$ and $\varphi$, we know that
$$
H(\nu): H_n(L_{\bullet}) \iso H_n(P_{\bullet}),
\quad 
H(\varphi): H_n(P_{\bullet}) \iso H_n(P_{\bullet}^{\#}) 
$$ 
and hence, the commutative diagram
$$\xymatrixcolsep{5pc}\xymatrix{
H_n(L_{\bullet}) \ar[r]^{H_n(\nu)}_{\sim} \ar[d]^{H_n(\mu_0)} & H_n(P_{\bullet})
\ar[d]^{H_n(\varphi)}_{\wr} \\
H_n(N_{\bullet}^{\#}) \ar[r]^{H_n(\nu_2^{\#})}_{\sim} & H_n(P_{\bullet}^{\#})
}$$ 
gives us (\ref{mu0iso}). We prove (\ref{otherH0}) now. Since the only nonzero
homology of $L_{\bullet}^{\#}$ is at degree $-n$, it is clear from the long exact
homology sequence for the bottom exact triangle that 
$$H_{i}(N_{\bullet}^{\#}) \cong H_{i}(P_{\bullet}^{\#}) ~\forall~i \neq -n, -n-1.$$
By (\ref{coneOfsubL}),
$H_i(N_{\bullet}) = 0$ for all $i\geq n$ and so
$$ H_{-(n+1)}(N_{\bullet}^{\#})= Ext^{d+2}\left(\frac{N_{n-1}}{B_{n-1}}, A \right) = 0,
 ~~  H_{-n}(N_{\bullet}^{\#})= Ext^{d+1}\left(\frac{N_{n-1}}{B_{n-1}}, A \right) = 0.$$
where $B_{n-1}\subseteq N_{n-1}$ is the boundary submodule
(the last part also follows directly because $Ext^d(\frac{P_n}{B_n}, A) \cong 
Ext^d(H_n(P_{\bullet}) , A)$). So, (\ref{otherH0}) is established. 
\pic $\eop$ 

\vspace{5mm} 
Now we consider the cone of $\mu_0$.
\begin{lemma}\label{RinInFPDcat}
With the notations in (\ref{leftSubLag}), (\ref{coneOfsubL}) and (\ref{homoOfM}),
consider an exact triangle on $\mu_0$ as follows:
$$
\Dia
L_{\bullet} \ar[r]^{\mu_0} & N_{\bullet}^{\#} \ar[r]^{\mu_1} 
& R_{\bullet} \ar[r]^{\mu_2} & T(L_{\bullet})
\enddiagram 
$$
where $R_{\bullet}$ is the cone of $\mu_0$. Then
$R_{\bullet}$ is  an object of $D^b_{{\mathcal A}}\left({\mathcal P}(A) \right)$.
More precisely,
$$
H_{i}(R_{\bullet}) =\left\{ 
\begin{array}{ll}
H_{i}(N_{\bullet}^{\#}) & for~~-(n-1) \leq i \leq n-1\\
0 & otherwise. 
\end{array}
\right. 
$$
which tells us that $R_{\bullet}$ has exactly two nonzero homologies less than
than $P_{\bullet}$.
\end{lemma} 
\pf 
Note that the only nonzero homology of $L_{\bullet}$ is at degree $n$. Using (\ref{homoOfM}),
the long exact homology sequence corresponding to the exact triangle is as follows: 

{\scalefont{0.7}
$$\xymatrixcolsep{1.3pc}
\xymatrix{
\cdots H_{n+2}(R_{\bullet})\ar[r] & 0\ar[r] & 0\ar[r] & H_{n+1}(R_{\bullet})\ar[r] & H_{n}(L_{\bullet})\ar[r]^{\sim} & H_{n}(N_{\bullet}^{\#})\ar[r] &  \\
H_{n}(R_{\bullet})\ar[r] & 0 \ar[r] & H_{n-1}(N_{\bullet}^{\#})\ar[r] & H_{n-1}(R_{\bullet})\ar[r]
& 0 \ar[r] & \cdots \cdots \ar[r] & \\ 
0 \ar[r] & H_{-n+1}(N_{\bullet}^{\#}) \ar[r] & H_{-n+1}(R_{\bullet}) \ar[r] & 0 \ar[r] & 0 \ar[r]
& H_{-n}(R_{\bullet}) \ar[r] & 0 \cdots
} 
$$}
Therefore, 
$$
H_{i}(R_{\bullet}) =\left\{ 
\begin{array}{ll}
H_{i}(N_{\bullet}^{\#}) & for~~-(n-1) \leq i \leq n-1 \\
0 & otherwise. 
\end{array}
\right. 
$$
\pic $\eop$ 
\begin{rem}\label{rmk3.5}
The readers are referred to \cite[4.11]{TWGI} for the definition of a
{\bf very good morphism} $L_{\bullet} \lra N_{\bullet}^{\#}$. All we need
in the sequel is that such morphisms exist \cite[4.17]{TWGI} and that
whenever they do, we obtain \cite[4.20]{TWGI}
\bE
\item There is a symmetric form $\psi:R_{\bullet} \iso
R^{\#}_{\bullet}$. 
\item There is a lagrangian 
$$
N_{\bullet}^{\#} \lra (P_{\bullet}^{\#}, \varphi^{-1}) \perp (R_{\bullet}, \psi).
$$ 
\eE 
\end{rem}
\bT\label{reduceL}
Let $(P_{\bullet}, \varphi)$ be a symmetric form as in (\ref{rmk3})
with $n > 0$. Then, there is a symmetric form $(Q_{\bullet}, \tau)$ such that
\bE 
\item 
$$
[(Q_{\bullet}, \tau)]=[(P_{\bullet}, \varphi)] \qquad in \quad 
W\left(D^b_{{\mathcal A}}\left(\P(A) \right)^{\pm}_u\right).
$$
\item $Q_{\bullet}$ has two less homologies  than 
$P_{\bullet}$ and it has  support in 
$[k+d, -k]$ for some $0\leq k < n$.
\eE 
\eT  
\pf Use the notations in (\ref{leftSubLag}),(\ref{coneOfsubL}) and (\ref{homoOfM}).
Using the above remark (\ref{rmk3.5}), let $\mu_0$ be a very good morphism and
let $(R_{\bullet}, \psi)$ be the symmetric form obtained.
Note that $N_{\bullet}^{\#}, R_{\bullet}$  are objects in 
$D^b_{{\mathcal M}FPD^{fl}}\left({\mathcal P}(A) \right)$  by (\ref{coneOfsubL})
and (\ref{RinInFPDcat}). Hence, by (2) of the above remark (\ref{rmk3.5}), 
$(R_{\bullet}, \psi)$ is Witt equivalent to $(P_{\bullet}^{\#}, -\varphi^{-1})$
in $\left(D^b_{{\mathcal A}}\left(\P(A) \right)^{\pm}_u\right)$ by definition (\ref{defWittGr}).
%$D^b_{{\mathcal M}FPD^{fl}}\left({\mathcal P}(A) \right)$. 
Therefore in 
$W\left(D^b_{{\mathcal A}}\left(\P(A) \right)^{\pm}_u\right)$,
we have
$$
[(R_{\bullet}, \psi)]=[(P{\bullet}^{\#}, -\varphi^{-1})] 
= [(P_{\bullet}, -\varphi)]  \quad 
\text{in} \quad 
W\left(D^b_{{\mathcal A}}\left(\P(A) \right)^{\pm}_u\right). 
$$
Now, $H_i(R_{\bullet}) = 0$ if $i < -(n-1), i > (n-1)$. By (\ref{stOfForms}),
$(R_{\bullet}, -\psi)$ is isometric to a form $(Q_{\bullet}, \tau)$ 
such that $Q_{\bullet}$ is supported on $[k+d, -k]$ with $0 \leq k \leq n-1$.
\pic $\eop$

\vspace{5mm}
Now we are ready to state and prove the main result of this article, which is
our version of the d\'{e}vissage theorem, i.e. 
\bT\label{IsoOfWzeta}
The homomorphisms
$$
W(\zeta): W\left({\mathcal A}, ^{\vee}, \pm\tilde{\varpi} \right)
\lra W\left(\left(D^b_{{\mathcal A}}\left(\P(A) \right)^{\pm}_u\right) \right)  
$$
%
%$$
%W(\zeta): W({\mathcal M}FPD^{fl}(A)) \lra  
%W^d\left(D^b_{{\mathcal M}FPD^{fl}}\left({\mathcal P}(A) \right)\right). 
%$$
induced by the functor 
$\zeta$
%: {\mathcal M}FPD^{fl}(A)   
%\lra D^b_{{\mathcal M}FPD^{fl}}\left({\mathcal P}(A) \right)$
are isomorphisms. 
\eT 
\pf First, we prove the surjectivity of the homomorphism $W(\zeta)$. Suppose
$x \in  W\left(\left(D^b_{{\mathcal A}}\left(\P(A) \right)^{\pm}_u\right) \right)$.
Then, by (\ref{stOfForms}), we can write $x=[(P_{\bullet},\varphi)]$ of the
form described in (\ref{rmk3}). Inductively, by (\ref{reduceL}), 
there is a form $(R_{\bullet},\psi)$ in 
$\left(D^b_{{\mathcal A}}\left(\P(A) \right)^{\pm}_u\right)$
such that 
$$
[(R_{\bullet}, \psi)]=[(P_{\bullet}, \varphi)]=x \qquad in \quad 
W\left(\left(D^b_{{\mathcal A}}\left(\P(A) \right)^{\pm}_u\right) \right)
$$
and $R_{\bullet}$ is supported in $[d,0]$. By (\ref{GoodAtLeft}), 
$R_{\bullet}$ is a projective resolution of  
$M:=H_0(R_{\bullet})\in {\mathcal A}$.
Further, $\psi$ induces a form $\psi_0:M \iso M^{\vee}$ and clearly
$$
W(\zeta)([(M,\psi_0)]) = [(R_{\bullet}, \psi)]=x.
$$
So, $W(\zeta)$ is surjective. 

Now, we proceed to prove that $W(\zeta)$ is injective. 
Suppose $(M,q)$ is a symmetric form in 
$\left({\mathcal A}, ^{\vee}, \pm\tilde{\varpi} \right)$ and 
$W(\zeta)([(M,q)])=0$. Write , $(\zeta(M), \zeta(q))=(P_{\bullet}, \varphi_0)$
where  $P_{\bullet}$ is a projective resolution of $M$ of length $d$,
and $\varphi_0$ is the induced symmetric form.
So, $[(P_{\bullet}, \varphi_0)]=0$ in 
$W\left(\left(D^b_{{\mathcal A}}\left(\P(A) \right)^{\pm}_u\right) \right)$.

This means there is a neutral form $[(Q_{\bullet}, \varphi_1)]$ so that
$[(P_{\bullet}, \varphi_0)] \perp [(Q_{\bullet}, \varphi_1)]$ is neutral
in $W\left(\left(D^b_{{\mathcal A}}\left(\P(A) \right)^{\pm}_u\right) \right)$.
Since $[(Q_{\bullet}, \varphi_1)]$ is neutral, so is $[(Q_{\bullet}, -\varphi_1)]$.
Hence, $[(P_{\bullet}, \varphi_0)] \perp [(Q_{\bullet}, \varphi_1)] \perp [(Q_{\bullet}, -\varphi_1)] $ is neutral. Using the usual isometry, we get that there is a hyperbolic form 
$$
\left(Q_{\bullet}\oplus Q_{\bullet}^{\#}, \left(\begin{array}{cc}
0 & 1 \\1 &0 \end{array} \right)\right)\qquad 
\text{with} \quad Q_{\bullet} \in 
D^b_{\A}(\left({\mathcal P}(A) \right)
$$
such that 
$$
(U_{\bullet}, \varphi):= 
\left(P_{\bullet} \oplus Q_{\bullet}\oplus Q_{\bullet}^{\#}, 
\left(\begin{array}{ccc}\varphi_0 & 0 & 0\\ 0& 0 & 1 \\
0& 1 &0 \end{array} \right)\right) 
\quad \text{is neutral in} \quad W\left(\left(D^b_{{\mathcal A}}\left(\P(A) \right)^{\pm}_u\right) \right)
$$
(we have so far followed the argument in \cite[3.5]{TWGI})\\
So, $(U_{\bullet}, \varphi)$ has a lagrangian $(L_{\bullet}, \alpha)$.
We have an exact triangle
$$
\diagram
T^{-1}L_{\bullet}^{\#} \ar[r]^w & L_{\bullet}  \ar[r]^{\alpha}  
& U_{\bullet}  \ar[r]^{\alpha^{\#}\varphi} 
& L_{\bullet}^{\#}  
\enddiagram 
\qquad T^{-1}w^{\#}=w.  
$$
Before proceeding, we use (\ref{easyHLofDual}) to make a startlingly simple observation
about the homologies of $L_{\bullet}^{\#}$ :
$$H_i(L_{\bullet}^{\#}) \cong H_{i-d}(L_{\bullet}^*) \cong Ext^d\left(H_{(d-i)-d}(L_{\bullet}) , A\right) \cong Ext^d\left(H_{-i}(L_{\bullet}) , A\right) \cong H_{-i}(L_{\bullet})^{\vee}.$$
Similarly, $H_i(U_{\bullet}^{\#}) \cong H_{-i}(U_{\bullet})^{\vee}$ and further using the
remark (\ref{rmk1.5nat}) about naturality of the maps, we get that the long exact homology sequences are
{\scalefont{0.7}
$$
\xymatrixcolsep{4pc}\xymatrix{
\ar[r] & H_{-2}(L_{\bullet})^{\vee} \ar[r]^{H_1(w)} \ar@{=}[d]
& H_{1}(L_{\bullet}) \ar[r]^{H_1(\alpha)} \ar[d]^{\wr}_{\tilde{\varpi}_{H_1}}
& H_1(U_{\bullet}) \ar[r]^{H_{-1}(\alpha)^{\vee} \smallcirc{0.5} H_1(\varphi)} 
\ar[d]^{\wr}_{H_1(\varphi)}
& H_{-1}(L_{\bullet})^{\vee} \ar@{=}[d] \cdots \cdots
& & &  \\
\ar[r] & H_{-2}(L_{\bullet})^{\vee} \ar[r]_{H_{-2}(w)^{\vee}} & H_1(L_{\bullet})^{\vee \vee}
 \ar[r]_{H_{-1}(\varphi)^{\vee} \smallcirc{0.5} H_1(\alpha)^{\vee \vee}} & H_{-1}(U_{\bullet})^{\vee} \ar[r]_{H_{-1}(\alpha)^{\vee}} & H_{-1}(L_{\bullet})^{\vee}
\cdots \cdots & & &  \\
}$$
$$\xymatrixcolsep{4pc}\xymatrix{
\cdots \cdots \ar[r]^{H_0(w)} &  H_0(L_{\bullet}) \ar[r]^{H_0(\alpha)} \ar[d]^{\wr}_{\tilde{\varpi}_{H_0}} &
H_0(U_{\bullet}) \ar[r]^{H_{0}(\alpha)^{\vee} \smallcirc{0.5} H_0(\varphi)} 
\ar[d]^{\wr}_{H_0(\varphi)}
& H_0(L_{\bullet})^{\vee} \ar[r]^{H_{-1}(w)} \ar@{=}[d] & \cdots \cdots  \\
\cdots \cdots\ar[r]_{H_{-1}(w)^{\vee}} & H_0(L_{\bullet})^{\vee \vee}
\ar[r]_{H_0(\varphi)^{\vee} \smallcirc{0.5} H_{0}(\alpha)^{\vee \vee}} &
H_0(U_{\bullet})^{\vee} \ar[r]_{H_0(\alpha)^{\vee}} &
H_0(L_{\bullet})^{\vee} \ar[r]_{H_0(w)^{\vee}} & \cdots \cdots  \\
}$$
$$\xymatrixcolsep{4pc}\xymatrix{
\cdots \cdots H_{-1}(L_{\bullet}) \ar[r]^{H_{-1}(\alpha)} \ar[d]^{\wr}_{\tilde{\varpi}_{H_{-1}}}
& H_{-1}(U_{\bullet}) \ar[r]^{H_1(\alpha)^{\vee} \smallcirc{0.5} H_{-1}(\varphi)} 
\ar[d]^{\wr}_{H_{-1}(\varphi)}
& H_{-1}(L_{\bullet})^{\vee} \ar[r]^{H_{-2}(w)} \ar@{=}[d] &
H_{-2}(L_{\bullet}) \ar[r]^{H_{-2}(\alpha)} \ar[d]^{\wr}_{\tilde{\varpi}_{H_{-2}}} & \\
\cdots \cdots H_{-1}(L_{\bullet})^{\vee \vee} \ar[r]_{H_1(\varphi)^{\vee} \smallcirc{0.5} 
H_{-1}(\alpha)^{\vee \vee}} &
H_1(U_{\bullet})^{\vee} \ar[r]_{H_1(\alpha)^{\vee}} &
H_1(L_{\bullet})^{\vee} \ar[r]_{H_1(w)^{\vee}} &
H_{-2}(L_{\bullet})^{\vee \vee} \ar[r]_{H_2(\varphi)^{\vee} \smallcirc{0.5} 
H_{-2}(\alpha)^{\vee \vee}} & \\
}
$$}
Replacing the part of the top exact sequence in negative degree by the corresponding
part of the bottom (dual) exact sequence, we get an exact sequence :
{\scalefont{0.7}
$$\xymatrixcolsep{5pc}\xymatrix{
\ar[r] & H_{-2}(L_{\bullet})^{\vee} \ar[r]^{H_1(w)} 
& H_{1}(L_{\bullet}) \ar[r]^{H_1(\alpha)} 
& H_1(U_{\bullet}) \ar[r]^{H_{-1}(\alpha)^{\vee} \smallcirc{0.5} H_1(\varphi)} 
& H_{-1}(L_{\bullet})^{\vee} \\
}$$
$$\xymatrixcolsep{5pc}\xymatrix{
\ar[r]^{H_0(w)} &  H_0(L_{\bullet}) \ar[r]^{H_0(\alpha)} &
H_0(U_{\bullet}) \ar[r]_{H_{0}(\alpha)^{\vee} \smallcirc{0.5} H_0(\varphi)} &
H_0(L_{\bullet})^{\vee} \ar[r]_{H_0(w)^{\vee}} & H_{-1}(L_{\bullet})^{\vee \vee} \\
}$$
$$\xymatrixcolsep{5pc}\xymatrix{
\ar[r]_{H_1(\varphi)^{\vee} \smallcirc{0.5} H_{-1}(\alpha)^{\vee \vee}} &
H_1(U_{\bullet})^{\vee} \ar[r]_{H_1(\alpha)^{\vee}} &
H_1(L_{\bullet})^{\vee} \ar[r]_{H_1(w)^{\vee}} &
H_{-2}(L_{\bullet})^{\vee \vee} \ar[r]_{H_2(\varphi)^{\vee} \smallcirc{0.5} 
H_{-2}(\alpha)^{\vee \vee}} & \\
}
$$
}
Notice that the complex above is very special, and is "symmetric" about $H_0(U_{\bullet})$.
So we can apply \cite[4.1 Lemma]{TWGII} to this sequence. Since the sequence is 
exact, we have
$$
[(H_0(U_{\bullet}), H_0(\varphi))]=[(0,0)]= 0
\qquad in \quad W({\mathcal M}FPD^{fl}(A)). 
$$ 
However,
$$
(H_0(U_{\bullet}), H_0(\varphi)) =
\left(H_0(P_{\bullet}) \oplus H_0(Q_{\bullet})\oplus H_0(Q_{\bullet}^{\#}), 
\left(\begin{array}{ccc}H_0(\varphi_0) & 0 & 0\\ 0& 0 & 1 \\
0& 1 &0 \end{array} \right)\right) 
$$
$$
= (M, q) \perp \left(H_0(Q_{\bullet})\oplus H_0(Q_{\bullet})^{\vee}, 
\left(\begin{array}{cc} 0 & 1 \\ 1 &0 \end{array} \right)\right)  
$$
So, we have 
$$
[(M, q)]= \left[(M, q) \perp \left(H_0(Q_{\bullet})\oplus H_0(Q_{\bullet})^{\vee}, 
\left(\begin{array}{cc} 0 & 1 \\ 1 &0 \end{array} \right)\right)
\right]
= 
[(H_0(U_{\bullet}), H_0(\varphi))]=0. 
$$
\pic $\eop$

%%%%%%%%%%%%%%%%%%%%%%%%%%%%%%%%%%%%%%%%%%%%%%%%%%%%%%%%%%%%%%%%%%%%%%%
%%%%%%%%%%%%%%%%%%%%%%%%%%%%%%%%%%%%%%%%%%%%%%%%%%%%%%%%%%%%%%%%%%%%%%%
\section{Shifted Witt Groups}\label{shiftedWSec} 
In this section, we use the previous results to obtain our d\'{e}vissage
theorem for the Witt groups $W^i(D^b_{\A}(\P(A)))$. We recall that $A$ is
a Cohen-Macauly ring with $\dim A_m = d \geq 2$ for all maximal ideals $m$
and such that $2$ is invertible in $A$ and that ${\mathcal A} = {\mathcal M}FPD^{fl}(A)$. 
\begin{notations}
For integers $j\geq 0$ define the functor
$\zeta_j = T^{-j} \smallcirc{0.7} \zeta$, which associates to an object
$M$ in ${\mathcal A}$ a projective resolution $P_{\bullet}$ of $M$ of length $d$,
such that $H_{-j}(P_{\bullet}) = M$.
\end{notations} 
\begin{defin}
 Suppose $K:=(K,\#,\delta,\varpi)$ is a triangulated category with
translation $T$ and $\delta-$duality $\#$.
We recall from \cite{TWGII} that 
$$T^n K := (K, T^n \smallcirc{0.5}\#,(-1)^n \delta),
(-1)^{\frac{n(n+1)}{2}} {\delta}^n {\varpi}).$$
is then also a triangulated category with the same translation $T$ but with
$((-1)^n \delta)$-duality $T^n \smallcirc{0.5} \#$.
If $K_0$ is a subcategory of $K$ satisfying the conditions of (\ref{defWittGr}),
we define $T^n K_0$ to be the same subcategory and translation with the induced
duality structure from $T^n K$. Using (\ref{defWittGr}), we define the shifted 
Witt groups by 
$$W^n(K):=W(T^nK)~~~~~~~~W^n(K_0):=W(T^nK_0)~~~~~~~\forall~ n \in {\mathbb Z}.$$
\end{defin}
Note that $T_s^2: T^n K \lra T^{n+4} K$ is 
an equivalence of triangulated categories with duality, for all 
$n\in{\mathbb Z}$. 
Similarly, $T_s^2: T^n K_0 \lra T^{n+4} K_0$ is an equivalence of categories
with duality, for all $n\in{\mathbb Z}$ and so
$$W^n(K) \iso W^{n+4}(K)~~~~~~~~W^n(K_0) \iso (T^nK_0)~~~~~~~\forall~ n \in {\mathbb Z}.$$
\begin{defin}
{\rm 
Following \cite{BW}, by {\bf "standard" duality} structure on ${\mathcal A}$,
we mean the exact categry 
$\left({\mathcal A}, ^{\vee}, (-1)^{\frac{d(d-1)}{2}}\tilde{\varpi}\right)$.  
By {\bf "standard" skew  duality} structure on ${\mathcal A}$,
we mean the exact categry
$\left({\mathcal A}, ^{\vee}, -(-1)^{\frac{d(d-1)}{2}}\tilde{\varpi}\right)$.
We denote the Witt groups  
$$
W_{St}^+(A)
=W\left({\mathcal A}, ^{\vee}, (-1)^{\frac{d(d-1)}{2}}\tilde{\varpi}\right), 
\quad
W_{St}^-(A)
=W\left({\mathcal A}, ^{\vee}, -(-1)^{\frac{d(d-1)}{2}}\tilde{\varpi}\right).  
$$ 
}
\end{defin} 
%%%%%%%%%%%%%%%%%%%%%%%%%%%%%%%%%%%%%%%%%%%%%%%%%%%%%%%%%%%%%%%%%%%%%%%
%%%%%%%%%%%%%%%%%%%%%%%%%%%%%%%%%%%%%%%%%%%%%%%%%%%%%%%%%%%%%%%%%%%%%%%
%%%%%%%%%%%%%%%%%%%%%%%%%%%%%%%%%%%%%%%%%%%%%%%%%%%%%%%%%%%%%%%%%%%%%%%
\bT\label{dshiftedWG}
Then, the functor $\zeta_0:{\mathcal A} \lra 
D^b_{{\mathcal A}}\left({\mathcal P}\left(A\right)\right)$ 
induces  an  isomorphism  
$$
W(\zeta_0): W_{St}^+(A) \iso 
W^d\left(D^b_{{\mathcal A}}({\mathcal P}(A)), *, 1, \varpi \right).
$$ 
\eT
\pf Recall $\zeta_0$ was denoted by $\zeta$ in the previous sections. 
For notational convenience
$\varpi_0=(-1)^{\frac{(d(d-1)}{2}}\varpi$.  
By theorem (\ref{IsoOfWzeta}), we get the following isomorphism of 
Witt groups 
$$
\eta_0:  W_{St}^+({\mathcal A}) \iso 
W\left(D^b_{{\mathcal A}}\left({\mathcal P}(A)\right), \#^u_d, 1, 
%(-1)^{\frac{d(d-1)}{2}}
\varpi_0  \right).  
$$
There is a duality preserving equivalence \cite[Proof of Lemma 6.4]{BW}
$$
\beta: \left(D^b\left({\mathcal P}(A)\right), T_u^d \smallcirc{0.5} *, 1, 
%(-1)^{\frac{d(d-1)}{2}}
\varpi_0  \right)
\lra 
\left(D^b\left({\mathcal P}(A)\right), T_u^d \smallcirc{0.5} * ,(-1)^d, (-1)^{\frac{d(d+1)}{2}}\varpi  \right).
$$
%Used $d+\frac{d-1}{2}=frac{d(d+1)}{2}$ 
Note that the later  is the shifted category 
$T^d\left(D^b\left({\mathcal P}(A)\right), * , 1, \varpi  \right)$.
Composing $\eta_0$ with the homomorphism induced by $\beta$, we get the result. $\eop$ 

%%%%%%%%%%%%%%%%%%%%%%%%%%%%%%%%%%%%%%%%%%%%%%%%%%%%%%%%%%%%%%%%%%%%%%%
%%%%%%%%%%%%%%%%%%%%%%%%%%%%%%%%%%%%%%%%%%%%%%%%%%%%%%%%%%%%%%%%%%%%%%%
%%%%%%%%%%%%%%%%%%%%%%%%%%%%%%%%%%%%%%%%%%%%%%%%%%%%%%%%%%%%%%%%%%%%%%%
%%%%%%%%%%%%%%%%%%%%%%%%%%%%%%%%%%%%%%%%%%%%%%%%%%%%%%%%%%%%%%%%%%%%%%%
\vspace{5mm} 
Now we prove the standard skew duality version of theorem (\ref{dshiftedWG}).
%%%%%%%%%%%%%%%%%%%%%%%%%%%%%%%%%%%%%%%%%%%%%%%%%%%%%%%%%%%%%%%%%%%%%%
\bT\label{skewWdm2}
The functor $\zeta_1$ induces an isomorphism
$$
W_{St}^-(A) \iso W^{d-2}\left(D^b_{{\mathcal A}}({\mathcal P}(A)), *, 1, -\varpi \right).  
$$
\eT 
\pf By Theorem (\ref{IsoOfWzeta}), we have an isomophism
$$
W_{St}^-(A) \iso W\left(D^b_{\mathcal A}({\mathcal P}(A)), T_u^d \smallcirc{0.5} *, 1 , -(-1)^{\frac{d(d-1)}{2}}\varpi   \right). 
$$
Write $\varpi_0= -(-1)^{\frac{d(d-1)}{2}}\varpi$. 
There is an equivalence of categories \cite[2.14]{TWGI}
$$
T_s:\left(D^b_{\mathcal A}({\mathcal P}(A)), T_u^{d-2} \smallcirc{0.5} *
, 1 ,%-(-1)^{\frac{d(d-1)}{2}}
\varpi_0 \right)   \lra 
\left(D^b_{\mathcal A}({\mathcal P}(A)), T_u^d \smallcirc{0.5} *, 1 ,
%-(-1)^{\frac{d(d-1)}{2}}
\varpi_0 \right).   
$$
This induces 
 an isomorphism 
$$
W\left(D^b_{\mathcal A}({\mathcal P}(A)), T_u^{d-2} \smallcirc{0.5} *, 1 ,
%-(-1)^{\frac{d(d-1)}{2}}
\varpi_0 \right)   \iso 
W\left(D^b_{\mathcal A}({\mathcal P}(A)), T_u^d \smallcirc{0.5} *, 1 ,
%-(-1)^{\frac{d(d-1)}{2}}
\varpi_0 \right).     
$$ 
%%%%%%%%%%%%%%%%
%By replacing $-\tilde{\varpi}$ by $-(-1)^{\frac{d(d-1)}{2}}\tilde{\varpi}$,
%in the proof of theorem \ref{Wzeta1WellD}, we have an isomorphims
%$$
%W(\zeta_{1}): W_s^-\left({\mathcal A}  \right) \iso
%W\left(D^b_{{\mathcal A}}
%\left({\mathcal P}(A), T_u, \#_{d-2}^u,1, -(-1)^{\frac{d(d-1)}{2}}\varpi  \right) \right).
%$$
%%%%%%%%%%%%%%%%%%%%%
As in the proof of \cite[Lemma 6.4]{BW}, there is an equivalence 
of triangulated categories with duality 
$$
\left(D^b({\mathcal P}(A)), T_u^{d-2} \smallcirc{0.5} *,1, 
%-(-1)^{\frac{d(d-1)}{2}}
\varpi_0  \right)  
\lra 
\left( D^b({\mathcal P}(A)), T_s^{d-2} \smallcirc{0.5} *,(-1)^{d-2}, 
(-1)^{\frac{d(d+1)}{2}}\varpi  \right). 
$$
The latter category is $T^{d-2} \left( D^b({\P}(A)), *, 1, -\varpi \right)$. 
\pic $\eop$ 
%%%%%%%%%%%%%%%%%%%%
\vspace{5mm}

Finally, we have the following regarding odd shifts.
\bT\label{oddShift}
For $n=d-1, d-3$, we have  
$$
W^{n}\left(D^b_{{\mathcal A}}({\mathcal P}(A)), *,1, \pm \varpi\right)=0.
%W^{d-3}\left(D^b_{{\mathcal A}}\left({\mathcal P}(A), T_s, ^*,1, -
%\varpi\right) \right)=0. 
$$
\eT 
\pf First consider $n=d-1$.
It would be enough to prove that 
$$
W\left(D^b_{{\mathcal A}}({\mathcal P}(A)), T_u^{d-1} \smallcirc{0.5} *, 1,
\pm \varpi\right) = 0.
%W\left(D^b_{{\mathcal A}}\left({\mathcal P}(A), T_u, \#_u^{d-3},1,\pm 
%\varpi\right) \right)=0. 
$$
Suppose $(P_{\bullet}, \varphi)$ is a form in 
$D^b_{{\mathcal A}}\left({\mathcal P}(A), T_u^{d-1} \smallcirc{0.5} *, 1, \pm \varpi\right)$.
By a little tweak in (\ref{stOfForms}), we can assume that 
$P_{\bullet}$ is supported on $[n+(d-1),-n]$ with $n > 0$
and $H_{-n}(P_{\bullet})\neq 0$. By imitating the arguments of theorem (\ref{reduceL}),
we can keep shortening the length of the complexes which give our symmetric form.
Eventually, we will be reduced to the case where the complex is 
$P_{\bullet}$ is supported on $[d-1, 0]$. By theorem (\ref{GoodAtLeft}),
$P_{\bullet}$ is exact. So, $[(P, \varphi)]=0$.
The same arguments apply when $n=d-3$. \pic $\eop$ 
%Now suppose  $n=d-3$ and suppose $(P, \varphi)$ is a form in 
%$D^b_{{\mathcal A}}\left({\mathcal P}(A), T_u,\#^u_{d-3},  ^*,1,
%\pm 
%\varpi\right)$. We can 
%assume that
%$P_{\bullet}$ is on $[n+(d-3),-n]$ with $n>0$
%and $H_{-n}(P_{\bullet})\neq 0$. We will prove $[(P, \varphi)]=0$.
%Again, by 
%arguments of
%theorem \ref{reduceL}, we can assume that
%$P_{\bullet}$ is supported on $[d-1, -2]$.
%By further application
%of the same argument, we can assume that
%$P_{\bullet}$ is supported on $[d-2, -1]$. By theorem \ref{GoodAtLeft},
%$P_{\bullet}$ is exact. So, $[(P, \varphi)]=0$.

\vspace{5mm} 
Using the $4$-periodicity, we now obtain the theorem mentioned in the introduction :
\bT[shiftFinal]\label{shiftFinal}
Let ${\mathcal B}=\left(D^b_{{\mathcal A}}({\mathcal P}(A)), T_s, ^*,1, 
\varpi\right)$. Then, for $n\in{\mathbb Z}$, we have   
\bE 
\item $W^{d+4n}\left({\mathcal B}\right)= W^+_s({\mathcal A})$,
\item $W^{d+4n+1}\left({\mathcal B}\right)= 0$, 
\item $W^{d+4n+2}\left({\mathcal B}\right)= W^-_s({\mathcal A})$,    
\item $W^{d+4n+3}\left({\mathcal B}\right)= 0$.    
\eE 
\eT 

\appendix 
\section{Some Formalism}\label{appA}
The purpose of this section is to prove the following theorem :
\bT\label{midWIsoF} Suppose $\E$ is a full subcategory of a $\Z [\frac{1}{2}]$ abelian category $\B$
with the $2$ out of $3$ property for short exact sequences, and has duality
$\left({\mathcal E}, ^{\vee}, \tilde{\varpi}\right)$. 
Let $D^b({\mathcal E}) := \left(D^b({\mathcal E}), *, a, \varpi \right)$ 
denote the derived category, with duality, of 
$\left({\mathcal E}, ^{\vee}, \tilde{\varpi}\right)$. Also, let 
$D^b_{{\mathcal E}}({\mathcal E})$ denote the derived category,
with duality, of objects in $D^b({\mathcal E})$ with homologies in 
${\mathcal E}$. Then the  homomorphism
$$
W(\mu): W\left({\mathcal E}, ^{\vee}, \tilde{\varpi} \right) \lra 
W\left(D^b_{{\mathcal E}}\left({\mathcal E} \right)  \right) 
$$
 induced by the functor    
$\mu: {\mathcal E} \lra D^b_{{\mathcal E}}\left({\mathcal E}  \right)$
is an isomorphism. 
\eT

In particular, with $\E = \A$ and $\B = {\mathcal M}(A)$, we obtain that
$$
W(\mu): W\left({\mathcal A}, ^{\vee}, \pm \tilde{\varpi} \right) \lra
W\left(D^b_{{\mathcal A}}
\left({\mathcal A}, ^{\vee}, \pm \tilde{\varpi} \right)\right) 
%W\left(D^b_{{\mathcal M}FPD^{fl}}
%\left({\mathcal M}FPD^{fl}(A) \right)  \right) 
$$
as promised in section \ref{WDIAGRAM}.

The proof of the theorem is essentially the same as the proof of
\cite[Theorem 3.2]{TWGII} with the extra check that all constructions
yield complexes whose homologies are in $\A$. This boils down to using
the most elementary of sublagrangians (concentrated in just one degree)
and reducing length. In any case, we follow the proof in \cite[Theorem 3.2]{TWGII}.
Since the category $\E$ has all the properties required in the results in 
\cite[Section 3]{TWGII}), we will freely borrow them.

To start with, injectivity of $W(\mu)$ follows directly because the isomorphism
$W\left({\mathcal E} \right) \iso W\left(D^b\left({\mathcal E} \right)  \right)$
(proven in  \cite[Theorem 4.3]{TWGII}) factors as
$$
\xymatrixcolsep{5pc}\xymatrix{
W\left({\mathcal E} \right) \ar[r]^{W(\mu)} \ar[rd]^{\sim} & 
W\left(D^b_{{\mathcal E}}\left({\mathcal E} \right) \right) \ar[d]\ar[d] \\
& W\left(D^b\left({\mathcal E} \right)  \right) 
}
$$ 
We move to the proof of surjectivity which, as we mentioned above will require
checking that we remain in $D^b_{\E}\left({\E} \right)$ through all the lemmas
establishing \cite[Theorem 3.2]{TWGII}. To start with, we establish the following
result regarding duality, which also provides an alternative proof of (\ref{relativeDUALITY}).
\begin{lemma}\label{serreDual}
With the same notations as in (\ref{midWIsoF}),
$D^b_{{\mathcal E}}\left({\mathcal E}\right)$ is closed under 
duality. 
\end{lemma}
\pf Let  
$P_{\bullet}$ be an object in the derived category
$D^b_{{\mathcal E}}\left({\mathcal E}\right)$.
Write  
$$
\diagram
P_{\bullet}:&\cdots \ar[r] &  P_2\ar[r]^{d_2} & P_1 \ar[r]^{d_1} & P_0 \ar[r]  & P_{-1} 
\ar[r] & P_{-2} \ar[r] & \cdots \\
P_{\bullet}^{*}:&\cdots\ar[r] &  P_{-2}^{\vee} \ar[r] &P_{-1}^{\vee}  \ar[r] &  P_0^{\vee}\ar[r] 
& P_1^{\vee} \ar[r] & P_2^{\vee} \ar[r]  & \cdots  
\enddiagram 
$$
Since the complexes are bounded and the homologies are objects of $\E$,
all the kernels $Z_i=ker(d_i)$, images $B_i=image(d_{i+1})$ and quotients 
$\frac{P_i}{B_i}$ are also objects of ${\mathcal E}$. Hence, so are their duals.
But we have an exact sequence 
$$
\diagram
0\ar[r] & \left(\frac{P_{t-1}}{B_{t-1}}\right)^{\vee}  \ar[r] 
&P_{t-1}^{\vee}  \ar[r] 
& \left(\frac{P_t}{B_t}\right)^{\vee} \ar[r] & H_t(P_{\bullet}^{*}) \ar[r]  & 0
\enddiagram 
$$
The first three terms in this sequence are in ${\mathcal E}$, hence so
is $H_0(P_{\bullet}^{*})$. \pic $\eop$ 
%%%%%%%%%%%%%%%%%%%%%%%%%%%%%%%%%%%%%%%%%%%%%%%%%%%%%%%%%%%%%%%%%%%%%%%
\begin{lemma} \label{3p6F} Let 
$x\in W\left(D^b_{{\mathcal E}}\left({\mathcal E}
\right)\right)$. Then $x=(P_{\bullet}, s)$ such that 
\bE
\item $P_{\bullet}$ is bounded and $s:P_{\bullet} \lra P_{\bullet}^{*}$ is a morphism
of complexes, without denominator.  
\item $s$ is quasi-isomorphism.
\item $s$ is strongly symmetric (i.e $s^{\vee}_{-i}=s_i~\forall i \in {\mathbb Z}$). 
\item $H_i(P_{\bullet}) \in {\mathcal E}$ for all 
$i \in {\mathbb Z}$. 
\eE
\end{lemma}
\pf
Let the form $x$ be given by $(X_{\bullet}, \eta)$. By definition there is a complex $P_{\bullet}$ which is an object of $D^b_{\E}(\E)$ and a chain complex
quasi-isomorphisms $t:P_{\bullet}\lra X_{\bullet}$ and
$\varphi_0: P_{\bullet}\lra X^*_{\bullet}$ such that $\eta=\varphi_0 t^{-1}$.
Then, $s=t^* \varphi_0 = t^* \eta t$ is a symmetric form on $P_{\bullet}$
and $(X_{\bullet}, \eta)$ is isometric to $(P_{\bullet}, \varphi)$. Clearly
$s$ is an actual morphism of complexes, a quasi-isomorphism and 
$H_i(P_{\bullet}) \in {\mathcal E}$ for all $i \in {\mathbb Z}$. Finally,
using that $\frac{1}{2}$ exists, we can make the map strongly symmetric.
$\eop$

\begin{lemma}\label{3p7F} 
Let $(P_{\bullet},s)$ be a symmetric form in $D^b_{\E}(\E)$ as in (\ref{3p6F}),
such that $P_{\bullet}$ is supported on $[m, -n]$ with $m > n \geq 0$.
Then $(P_{\bullet}, s)$ is isometric to a symmetric space $(Q_{\bullet}, t)$
such that $Q_{\bullet}$ is supported on $[n,-n]$ and $(Q_{\bullet}, t)$ has
all the other properties of $(P_{\bullet},s)$.
\end{lemma}
\vspace{5mm} 
\pf This is precisely \cite[Lemma 3.7]{TWGII} in our context. 
Note that since we can do this in the derived category without the homology
condition $D^b(\E)$, we use the same result to get $(Q_{\bullet}, t)$ isometric
to $(P_{\bullet},s)$ with the required condition. However, since the isometry
gives in particular a quasi-isomorphism $P_{\bullet} \iso Q_{\bullet}$ and
so $Q_{\bullet}$ is also an object in $D^b_{\E}(\E)$.
\pic $\eop$  

%%%%%%%%%%%%%%%%%%%%%%%%%%%%%%%%%%%%%%%%%%%%%%%%%%%%%%%%%%%%%%%%%%%%%%% 
\begin{lemma} \label{3p9F}
Let $(P_{\bullet}, s)$ be a symmetric space, as in (\ref{3p6F}).
with support on $[-n, n]$ and $n > 0$. Then there exists
a symmetric space $(Q_{\bullet},t)$ such that
\bE 
\item $(Q_{\bullet},t)$ is as in (\ref{3p6F}).
\item $(Q_{\bullet},t)$ is supported in $[n, -(n-1)]$. 
\item $H_i(Q_{\bullet}) \in {\mathcal E}$
for all $i \in {\mathbb Z}$.
\item $[(P_{\bullet},s)]+[(Q_{\bullet},t)]=0$ in 
$W\left(D^b_{{\mathcal E}}\left({\mathcal E}\right)\right)$.
\eE
\end{lemma} 
\pf 
Once again this is \cite[Lemma 3.9]{TWGII} in our context.
We skim through the proof mentioning only the significant points and most
important, checking the points where we need to check the extra homology condition.
We begin by proving the lemma in the case $n \geq 2$. Write 
$$
\diagram
P_{\bullet}=\ar[d]_s &
\cdots~0 \ar[r] &  P_{n} \ar[d]^s\ar[r]^d&P_{n-1}\ar[r]^d\ar[d]_s &\cdots \ar[r]^d 
& P_{-n} \ar[d]^s\ar[r]  &0 \\
P_{\bullet}^{*}=& 
\cdots~0 \ar[r] &  P_{-n}^{\vee} \ar[r]_{d^{\vee}}
&P_{-(n-1)}^{\vee}\ar[r]_{d^{\vee}}&\cdots \ar[r] &   P_{n}^{\vee}  \ar[r] & 0  \\ 
\enddiagram 
$$
Define $(Q_{\bullet}, t)$ as follows, on the left side:  
{\scalefont{0.65}
$$
\xymatrixcolsep{1.3pc}\xymatrix{
Q_{\bullet}=\ar[d]_t 
\cdots~0 \ar[r] &  P_{n} \ar[d]^0\ar[r]^{\left(\begin{array}{c} s\\d\end{array} \right)\quad}
&P_{-n}^{\vee}\oplus P_{n-1}\ar[r]^{(0,d)}\ar[d]_{(d^{\vee},-s)} &P_{n-2}\ar[d]_{-s} \ar[r]^d & \cdots P_{-(n-2)}\ar[d]_{-s} \ar[r]^d & P_{-(n-1)}\ar[r]
\ar[d]^{\left(\begin{array}{c} d\\-s\end{array} \right)}
& 0 \ar[r]  &0 \\
Q_{\bullet}^{*}=
\cdots~0 \ar[r] &  0 \ar[r]_{d^{\vee}}
&P^{\vee}_{-(n-1)}\ar[r]_{d^{\vee}}&P^{\vee}_{-(n-2)}\ar[r]_{d^{\vee}}  &\cdots 
P^{\vee}_{(n-2)}\ar[r]_{\left(\begin{array}{c} 0\\d^{\vee}\end{array} \right)}  &P_{-n}\oplus P^{\vee}_{n-1} \ar[r]_{(s,d^{\vee})} &   P_{n}^{\vee}  \ar[r] & 0  \\ 
}$$
}
It was proved in \cite[Lemma 3.9]{TWGII} that $t$ is a quasi-isomorphism. 
So, it follows
$$
H_n(Q_{\bullet}) \cong 0, \quad H_{n-1}(Q_{\bullet}) \cong H_{n-1}(Q_{\bullet}^{*})
\cong \ker(d^{\vee}) \in {\mathcal E}.  
$$
Since $\text{image}(0,d_{n-1}) = \text{image}(d_{n-1})$ we have
$$
H_{i}(Q_{\bullet})=H_{i}(P_{\bullet}) \in {\mathcal E} 
\qquad \forall \quad i\leq n-2. 
$$
Therefore
$$
H_{i}(Q_{\bullet}) \in {\mathcal E}
\qquad 
for ~all \quad i \in {\mathbb Z}. 
$$
So $Q_{\bullet}$ satisfies the last condition of (\ref{3p6F}). 
The other conditions of (\ref{3p6F}) are shown to be established in \cite[Lemma 3.9]{TWGII}.

Therefore, $Q_{\bullet}$ satisfies (\ref{3p6F}).
It was established in \cite{TWGII} that 
$(P_{\bullet}, s) \perp (Q_{\bullet}, t)$ is neutral in 
$D^b\left({\mathcal E}\right)$, by showing that 
$(P_{\bullet}, s) \perp (Q_{\bullet}, t)$
is isometric (in $D^b\left({\mathcal E}\right)$) to the cone 
of the morphism $z:T^{-1}M_{\bullet}^{*} \lra M_{\bullet}$ 
defined as follows:
$$
\diagram
T^{-1}M_{\bullet}^{*}=\ar[d]_z
\cdots
0\ar[r] &  P_{n} \ar[d]^0\ar[r]^{-d}&P_{n-1}\ar[r]^{-d}\ar[d]^0 
&\cdots P_{-(n-2)}
\ar[d]^0\ar[r]^{-d} 
& P_{-(n-1)} \ar[d]^{s_{-n}d}\ar[r]  &0 \\
M_{\bullet}= 
\cdots 0 \ar[r]&  P_{-(n-1)}^{\vee} \ar[r]_{d^{\vee}}
&P_{-(n-2)}^{\vee}\ar[r]_{d^{\vee}}&\cdots  P_{n-1}^{\vee}\ar[r] &   P_{n}^{\vee}  \ar[r] & 0  \\ 
degree=& 
n\ar@{-->}[u] &  n-1 \ar@{-->}[u]  
&& &   -n  \ar@{-->}[u] &   \\ 
\enddiagram 
$$
Since all the boundaries and cycles of $P_{\bullet}$ and $P_{\bullet}^{*}$
are objects of $\E$, so are $H_i(M_{\bullet})$ and $H_i(M_{\bullet})^*$.
Therefore,$M_{\bullet}$ and $M_{\bullet}^*$ are objects of 
$D^b_{{\mathcal E}}\left({\mathcal E}\right)$.

Let $Z_{\bullet}=cone(z)$. In \cite[Lemma 3.9]{TWGII}, it is further shown
that there is a symmetric form $\chi:Z_{\bullet}\lra Z_{\bullet}^{*}$ and
that $(Z_{\bullet}, \chi)$ is isometric to $(P_{\bullet},s) \perp (Q_{\bullet},t)$
in $D^b\left({\mathcal E} \right)$. But this tells us that 
$H_i(Z_{\bullet}) \cong H_i(P_{\bullet}) \oplus H_i(Q_{\bullet})$ and hence
$Z_{\bullet}$ is an object of $D^b_{{\mathcal E}}\left({\mathcal E}\right)$.

Now again following the proof of \cite[Lemma 3.9]{TWGII}, it is shown that
$T^{-1}z^{\#} = z$ in $D^b_{{\mathcal E}}\left({\mathcal E} \right)$ and that
the form $\chi$ actually fits in to make $M_{\bullet}$ a lagrangian for
$(Z_{\bullet}, \chi)$. Hence, this proves the lemma when $n\geq 2$.

In the case $n=1$, $(P_{\bullet}, s)$ is given by 
$$
\diagram 
P_{\bullet}=\ar[d]_s&
0\ar[r] & P_1\ar[d]_s \ar[r] & P_0\ar[d]_s \ar[r] & P_{-1} \ar[r]\ar[d]_s & 0 \\
P_{\bullet}^{*}=&
0\ar[r] & P_{-1}^{\vee} \ar[r] & P_0^{\vee} \ar[r] & P_{1}^{\vee} \ar[r] & 0 \\
\enddiagram 
$$
Define $(Q_{\bullet}, s)$  as follows 
$$
\xymatrixcolsep{4pc}\xymatrixrowsep{4pc}\xymatrix{
Q_{\bullet} = \cdots 0 \ar[d]_t  \ar[r] & P_1\ar[d] 
\ar[r]^{\left(\begin{array}{c}s \\ d \end{array} \right)~~~} 
&P_{-1}^{\vee}\oplus  P_0\ar[d]^{\left(\begin{array}{cc}0 & d \\d^{\vee}&-s \end{array} \right)} \ar[r] & 0 \ar[r]\ar[d] & 0 \\
Q_{\bullet}^{*} = \cdots 0 \ar[r] 
 & 0 \ar[r] & P_{-1}\oplus P_0^{\vee} 
\ar[r]_{\left(\begin{array}{cc}s^{\vee} d^{\vee} \end{array} \right)} 
& P_{1}^{\vee} \ar[r] & 0 \\
}
$$
The degree zero term is in the middle. In \cite[Lemma 3.9]{TWGII}, it is established
that $t$ is a quasi-isomorphism. It follows that $H_i(Q_{\bullet})=0$ for all $i\neq 0$
and 
$$
H_0(Q_{\bullet})=\frac{P_{-1}^{\vee}\oplus  P_0}{P_1} \in {\mathcal E}.  
$$
So, $Q_{\bullet}$ satisfies all the condition in (\ref{3p6F}), because the
remaining three conditions are established in \cite[Lemma 3.9]{TWGII}.
Now $(P,s) \perp (Q, t)$ has a lagrangian, namely
$$
\diagram
T^{-1}M_{\bullet}^{*}=\ar[d]_z& 0 \ar[r] & P_1 \ar[r]^{-d}\ar[d]^0 
& P_{0}\ar[r]\ar[d]^{sd} & 0\\  
M_{\bullet}=& 0 \ar[r] & P_0^{\vee} \ar[r]_{d^{\vee}} & P_{1}^{\vee}\ar[r] & 0\\  
degree=& 1\ar@{-->}[u] & 0 \ar@{-->}[u] & -1 \ar@{-->}[u] & 
\enddiagram 
$$
Again,
$$
H_0(M_{\bullet}) =\ker(d^{\vee}),  \quad  
H_1(M_{\bullet}) =\text{coker}(d^{\vee})\quad are ~in \quad {\mathcal E}.
$$
So, $M_{\bullet}$ and hence $M_{\bullet^*}$ are objects of
$D^b_{{\mathcal E}}\left({\mathcal E} \right)$. The rest of the proof is
exactly the same as in the case $n \geq 2$.
\pic $\eop$ \\

{\bf \noindent Finishing the prof of (\ref{midWIsoF})} : \\
We use (\ref{3p6F}) to represent any element $x$ in $W(D^b_{\A}(\A))$
by a chain complex in $D^b_{\A}(\A)$ and a strongly symmetric quasi-isomorphism
to its dual. Then, by alternate use of lemma (\ref{3p7F}) and (\ref{3p9F}), 
we reduce any element in $W\left(D^b_{{\mathcal E}}\left({\mathcal E} \right)
\right)$ to a chain complex in $D^b_{{\mathcal E}}\left({\mathcal E}) \right)$,
concentrated at degree zero. Of course that means the quasi-iromorphism is actually
an isomorphism and hence $x$ is the image of an element in $W(\A)$ via $W(\mu)$.
So $W(\mu)$ is also surjective.
So, the proof of theorem (\ref{midWIsoF}) is complete. 
$\eop$
% \newpage 

\end{document}